\documentclass[12pt,sort&compress]{elsarticle}

\usepackage[pdfborder={0 0 0},colorlinks,allcolors=blue]{hyperref}

\usepackage{mystyle}
\usepackage{afterpage}

\theoremstyle{definition}

\journal{}

\begin{document}

\begin{frontmatter}

\title{Interpolation-based reproducing kernel particle method}

\author[ucsd]{Jennifer E. Fromm}
\author[colorado]{John A. Evans} \author[ucsd]{J. S. Chen\corref{cor1}}
\ead{js-chen@ucsd.edu}
\cortext[cor1]{Corresponding author}

\address[ucsd]{University of California San Diego, San Diego, CA 92093, USA}
\address[colorado]{University of Colorado Boulder, Boulder, CO 80309, USA}

\begin{abstract}
Meshfree methods, including the reproducing kernel particle method (RKPM), have been widely used within the computational mechanics community to model physical phenomena in materials undergoing large deformations or extreme topology changes. 
RKPM shape functions and their derivatives cannot be accurately integrated with the Gauss-quadrature methods widely employed for the finite element method (FEM) and typically require sophisticated nodal integration techniques, preventing them from easily being implemented in existing FEM software. 
Interpolation-based methods have been developed to address similar problems with isogeometric and immersed boundary methods, allowing these techniques to be implemented within open-source finite element software. 
With interpolation-based methods, background basis functions are represented as linear combinations of Lagrange polynomial foreground basis functions defined upon a boundary-conforming foreground mesh. 
This work extends the applications of interpolation-based methods to implement RKPM within open-source finite element software. 
Interpolation-based RKPM is applied to several PDEs, and error convergence rates are equivalent to classic RKPM integrated using high-order Gauss-quadrature schemes.
The interpolation-based method is able to exploit the continuity of the RKPM basis to solve higher-order PDEs, demonstrated through the biharmonic problem. 
The method is extended to multi-material problems through Heaviside enrichment schemes, using local foreground refinement to reduce geometric integration error and achieve high-order accuracy. 
The computational cost of interpolation-based RKPM is similar to the smoothed gradient nodal integration schemes, offering significant savings over Gauss-quadrature-based meshfree methods while enabling easy implementation within existing finite element software.. 

\end{abstract}


\end{frontmatter}

\section{Introduction}
\label{sec:introduction}

The finite element method (FEM) is a fundamental tool within computer aided engineering (CAE) and is widely used in both industrial design and within academic settings. Within FEM, a geometric domain is discretized into a mesh, which is also used to generate basis functions to approximate the solutions to partial differential equations (PDEs) \cite{hughes_finite_2000}. While well suited to many applications, FEM is limited by the requirements of meshing, which can consume up to 80\% of an engineering analyst's time \cite{boggs_dart_2005}. FEM simulation quality is also closely tied to mesh quality \cite{burkhart_finite_2013,knupp_remarks_2007}, especially when high-order methods are employed \cite{engvall_mesh_2020}.

Meshfree methods encompass a large group of computational tools that approximate PDEs with bases constructed on scattered points not connected through a mesh \cite{chen_meshfree_2017,belytschko_meshfree_2023}. Meshfree methods stemmed from the diffuse element method \cite{nayroles_generalizing_1992}, which was soon followed by the development of element free Galerkin (EFG) methods \cite{belytschko_element-free_1994,belytschko_crack_1995,lu_new_1994}. Drawing from wavelet methods, the reproducing kernel particle method (RKPM) was introduced in \cite{liu_reproducing_1995} and constructed with discrete reproducing conditions to ensure optimal convergence and extended to large deformation problems in \cite{chen_reproducing_1996}. In addition to circumventing meshing (and re-meshing, in the case of large deformation problems) challenges, meshfree methods offer several advantages over classical FEM. Meshfree methods in general, and RKPM in particular, allow for convenient local h-refinement through node insertion when stabilized conforming nodal integration (SCNI) \cite{chen_stabilized_2001,chen_non-linear_2002} was used. $p-$refinement, increasing the polynomial order of basis functions to improve accuracy is also possible with meshfree methods.  While $p-$refinement increases the cost of evaluating RKPM basis functions, it does not increase the number of system degrees of freedom (DOFs) as it does in the case of FEM. And perhaps most notably, the continuity of RKPM bases are independent of polynomial order and controlled solely by the kernel. 

Numerical integration of basis functions over domain geometry has historically presented a major challenge in the implementation of RKPM. 
The basis functions used in classic finite element methods are compactly supported and can be accurately integrated using Gauss quadrature rules within a desirable order of precision. Gauss quadrature has also be used for meshfree methods \cite{dolbow_numerical_1999, lu_new_1994}.  As some meshfree shape functions, including RKPM, are by nature rational functions with overlapping supports incompatible with Gauss quadrature integration cells, the number of Gauss quadrature points required to reduce integration error to below approximation error can result in impractical computation times. 

Additionally, evaluation algorithms require a mechanism to identify the nearest neighbors of each node to compute basis support size and determine non-zero shape functions at each evaluation point, increasing complexity. While these tasks can be performed off-line as preprocessing steps, reducing the number of evaluation points is critical to implementing efficient RKPM. 
Nodal integration schemes with stabilization have thus been developed to enable large scale computation of meshfree methods \cite{chen_stabilized_2001,chen_arbitrary_2013,puso_meshfree_2008,hillman_accelerated_2016}. 
As these integration tools are quite different from those used in conventional mesh-based methods, meshfree methods still require unique implementation software extensions and have not been widely adopted within commercial simulation codes.

Other classes of computational mechanics methods that were developed to avoid mesh generation face similar integration challenges. Isogeometric methods \cite{hughes_isogeometric_2005} utilize B-spline basis functions which are similarly rational and largely supported, making integration difficult.  Immersed boundary methods \cite{peskin_flow_1972,parvizian_finite_2007,schillinger_finite_2015}, including immersed isogeometric (often called immersogeometric) \cite{kamensky_immersogeometric_2015, noel_xiga_2022, schmidt_extended_2023} methods, lack a boundary conforming mesh to define quadrature rules upon. 

Interpolation has previously been applied to address similar challenges within isogeometric analysis \cite{schillinger_lagrange_2016,tirvaudey_non-invasive_2020,kamensky_tigar_2019,kamensky_open-source_2021}and immersed boundary methods \cite{fromm_interpolation-based_2023,fromm_interpolation-based_2024}. Interpolation-based methods utilize extraction to represent background basis functions that are difficult to integrate as linear combinations of foreground basis functions that are easy to integrate. For example, with B\'ezier extraction, the B-splines of an isogeometric basis are exactly represented as linear combinations of Bernstein polynomials \cite{scott_isogeometric_2011}. This was generalized to Lagrange extraction in \cite{schillinger_lagrange_2016}, which was then utilized to implement isogeometric analysis via classic finite element software as in \cite{kamensky_open-source_2021}. 

Implementations of B\'ezier and Lagrange extraction typically utilize exact interpolation, where a background B-spline basis function is exactly represented by a foreground basis. This places limitations on the construction of the foreground basis. For example, to exactly interpolate a bi-quadratic B-spline patch with Lagrange polynomials defined on triangular elements, the edges of the triangular elements must conform with the edges of the patch, and fourth-order polynomials must be employed. 

Interpolation-based immersed boundary methods introduced the concept of ``approximate extraction", where the interpolated basis is not everywhere equivalent to the target background basis \cite{fromm_interpolation-based_2023}. 
This approximation results in interpolation error. However, the error is limited to the regions surrounding material or domain interfaces and for exact geometric representation is bounded by method error. For this reason interpolation-based immersed boundary methods can still maintain ideal error convergence rates. 

Venturing beyond what has been mathematically shown, \cite{fromm_interpolation-based_2023} also presents numerical evidence that approximate extraction can be extended to loosen the constraints placed on the foreground mesh. Non-background conforming foreground meshes allowed the supports of background basis functions to arbitrarily intersect the foreground elements. This work builds upon the concept of approximate extraction and non-background fitted interpolation with the extension of interpolation-based methods to RKPM. This method will be referred to as Interpolation-based RKPM, or Int-RKPM. 

The outline of this paper is as follows: Section \ref{sec:RKPM} provides an overview of RKPM and existing integration techniques; Section \ref{sec:Int-RKPM} introduces the novel computational method combining interpolation with RKPM and demonstrates its efficacy when compared to classic RKPM; Section \ref{sec:localRef} expands the method with local foreground refinement strategies to model problems on curved geometries; Section \ref{sec:heaviside} implements the method with enrichment strategies to model multi-material problems; Finally Section \ref{sec:conclusion} concludes the work and indicates future directions and applications of Int-RKPM. 

\section{An overview of the reproducing kernel particle method}\label{sec:RKPM}

RKPM is a meshfree method for solving PDEs. Instead of creating a mesh data structure, a domain $\Omega \subset \mathbb{R}^d $ is discretized by a set of $NP$ nodes of spatial coordinates $ \{ \bm{x}_{1}, \bm{x}_{2}, ... \bm{x}_{NP}\}$. Here $d$ denotes the spatial dimension. For ease of explanation, this work will deal only with the case $d=2$, with notation $\bm{x} = [x, y]$. RKPM shape functions $\Psi_{I}$ are associated with each node $\bm{x}_I$.
The following subsections will detail the construction and evaluation of these shape functions and review integration techniques used in classic RKPM implementations. 

\subsection{RKPM shape functions}
 The RKPM approximation of a function $u$ is given as 
\begin{align}
    u(\bm{x}) \approx u^h (\bm{x})  = \sum_{I} \Psi_{I} (\bm{x}) u_I,
\end{align}
where each $\Psi_{I}$ is a shape function and $u_I$ can be collected into a vector of coefficients. The shape functions can be expressed as 
\begin{align}\label{eq:RKPMShapeCOmp}
    \Psi_{I}(\bm{x}) = \Phi_{a}(\bm{x} - \bm{x}_{I}) C(\bm{x}; \bm{x} -  \bm{x}_{I}), 
\end{align} 
where $\Phi_a (\bm{x} )$ is the kernel function with compact support size $a_{I}$ defined as 
\begin{align}
    a_{I} = c_a h_I, 
\end{align}
where $c_a$ is the normalized support size and $h_I$ is the nodal spacing associated with the node $\bm{x}_{I}$. The nodal spacing is determined by 
\begin{align}
    h_{I} = \text{max} \{ || \bm{x}_I - \bm{x}_J|| \}, \, \forall \bm{x}_{J} \ \in B_{I}
\end{align}
where $B_{I}$ is the set of closest nodes to $\bm{x}_I$. The number of nodes contained within $B_I$ is a tunable parameter based upon dimension and domain specificiations. For this study with $d=2$, $B_I$ contains four nodes. The kernel function $\Phi_a$ determines the local order of continuity. The kernel functions support shape also dictates the shape function's support shape, and is typically either circular or rectangular for $d=2$. 

The function $C(\bm{x}; \bm{x} - \bm{x}_{I})$ in Equation \ref{eq:RKPMShapeCOmp} is the correction function and is used to impose the reproducing conditions for the desired $\Bar{\Omega}$ completeness. The reproducing condition states that 
\begin{align}
    \sum_{I} \Psi_{I} (\bm{x}) x^i_I y^j_I =  x^i y^j, \qquad 0 \leq i + j \leq n
\end{align}
where $n$ is the order of reproducibility. The correction function can be computed as 
\begin{align}\label{eq:RKPMShapeCOmp}
    C(\bm{x}; \bm{x} - \bm{x}_{I})= \bm{H}^{\text{T}}(0) \bm{M}^{-1}(\bm{x})\bm{H}(\bm{x} - \bm{x}_{I}),
\end{align} 
where $\bm{H}(\bm{x})$ is the basis vector 
\begin{align}
    \bm{H}(\bm{x}) = [ 1, x, y, x^2, xy, y^2, ..., x^n, x^{n-1}y, ..., x y^{n-1},  y^n ], 
\end{align}
and $\bm{M}(\bm{x})$ is the moment matrix, 
\begin{align}\label{eq:momentMat}
    \bm{M}(\bm{x}) = \sum_{I \in G_{\bm{x}}} \bm{H}(\bm{x} - \bm{x}_I) \bm{H}^{\text{T}}(\bm{x} - \bm{x}_I) \Phi_{a}(\bm{x} - \bm{x}_{I}),
\end{align}
where $G_{\bm{x}}$ is the set of all nodes with support covering the point $\bm{x}$. To avoid singular moment matrices, the support of each shape function  (and thus kernel function) must be sufficiently large to provide adequate coverage to every evaluation point within the domain. Each point must be covered by at least $n_p$ kernels, where 
\begin{align}
    n_p = \begin{pmatrix}
        n + d \\
        d
    \end{pmatrix}. 
\end{align}
This is ensured by altering the normalized support size $c_a$. 

\begin{figure}
\begin{subfigure}[b]{0.5\linewidth}
	 \centering
    \includegraphics[width=\linewidth]{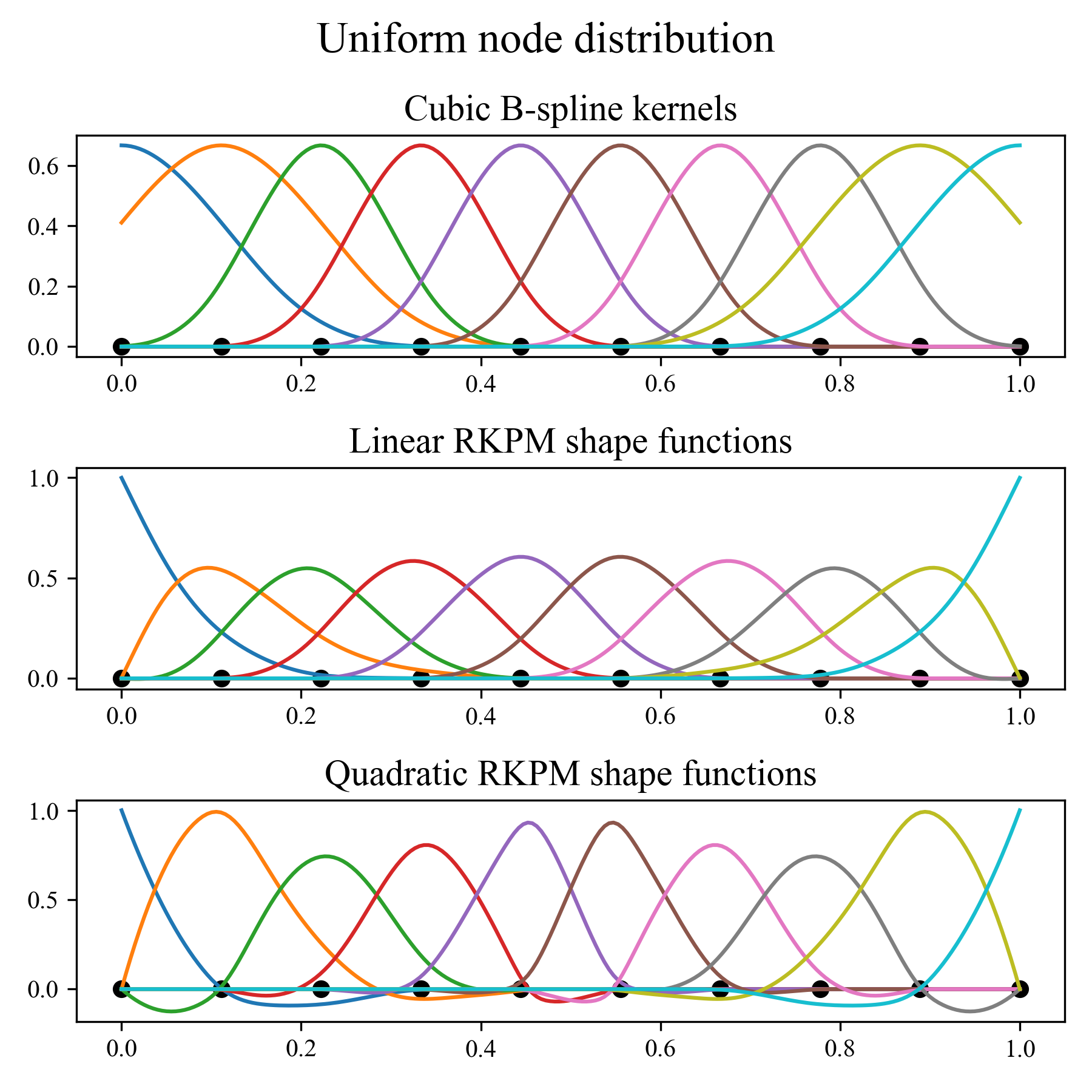}
    \caption{\label{fig:RKPM1D_uniform}1D uniform nodal distribution. }
  \end{subfigure}
	 \begin{subfigure}[b]{0.5\linewidth}
	 \centering
    \includegraphics[width=\linewidth]{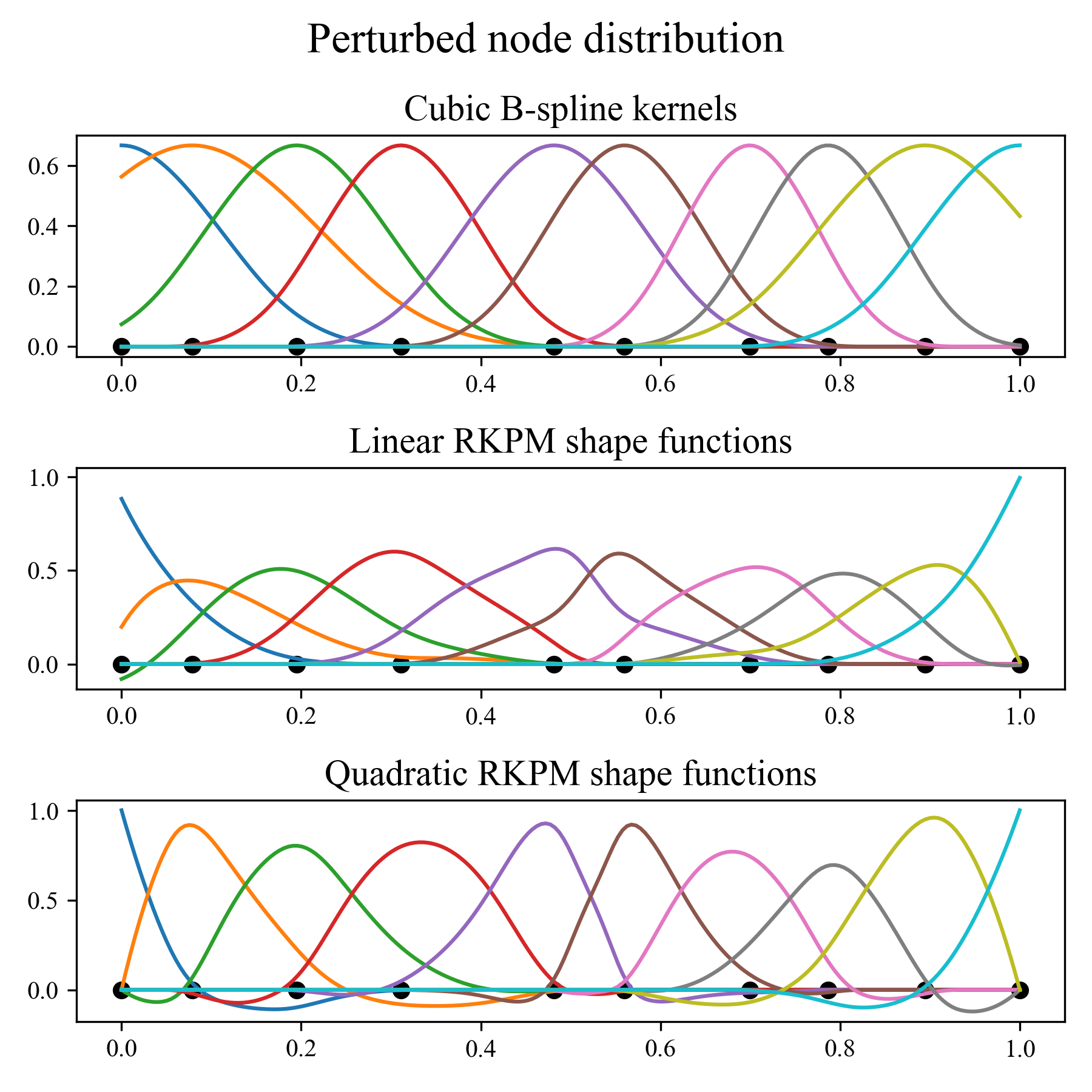}
    \caption{\label{fig:RKPM1D_rand} 1D non-uniform nodal distribution.}
  \end{subfigure}
  \caption{ Cubic B-spline kernels (top), linear RKPM shape functions (middle) and quadratic RKPM shape functions (bottom) are plotted. With a uniform nodal distribution (left), the shape functions are very regular. Even with a non-uniform distribution (right) the functions satisfy the reproducing condition. }
\end{figure}

Examples of $d=1$ RKPM shape functions are plotted in Figures \ref{fig:RKPM1D_uniform} and \ref{fig:RKPM1D_rand}, for both a uniform and non-uniform (perturbed) nodal distribution. Cubic B-spline kernels are employed, with support sizes such that each kernel covers at least three neighboring nodes. Both linear and quadratic reproducing conditions are applied, to compute the associated RKPM shape functions.

The evaluation of a given shape function at a point within the computational domain requires the inversion of the moment matrix $\bm{M}(\bm{x})$. While this matrix is small relative to the linear systems typically dealt with in computational mechanics, the inversion must be performed at every evaluation point, for every shape function supporting that point. Thus efficient assembly workflows that minimize the number of evaluation points are essential for the implementation of RKPM.

The product rule is used to compute the derivatives of shape functions 
\begin{align}\label{eq:RKPMShapeCOmp}
    \dfrac{\partial \Psi_{I}}{\partial x} = \dfrac{\partial \Phi_{a}}{\partial x} C+  \Phi_{a}  \dfrac{\partial C}{\partial x}
\end{align} 
where computing the derivative of the correction function requires additional inversions of the derivatives of moment matrices. This cost is further exacerbated for second and third derivatives, which are required for the approximation of higher than second order PDEs such as the biharmonic problem investigated in Subsection \ref{subsec:biharmonic}. 

\subsection{RKPM integration techniques }
Evaluating RKPM shape functions at a given point (Equation \ref{eq:RKPMShapeCOmp}) requires the inversion of the moment matrix, making each evaluation computationally expensive. This computational expense can be made up for by the relative reduction in over-all degrees of freedom when compared with classical FEM, but remains a motivating factor in the development of more efficient integration schemes. Gaussian integration, performed on structured integration cells can be utilized, but in addition to requiring large numbers of evaluation points this technique yields issues in accuracy and also negatively affect convergence. Nodal integration schemes have thus become the most popular method for implementing RKPM. 

Direct nodal integration (DNI), where shape functions and their derivatives are only evaluated at nodes presents challenges as the first derivatives of meshfree shape functions vanish at nodes (or nearly vanish at nonuniformly distributed nodes). Vanishing or nearly vanishing derivatives result in an oscillatory mode of wavelength $2h$ (where $h$ is the nodal spacing) leading to instability \cite{hillman_accelerated_2016, belytschko_meshfree_2023}. 
Numerous stabilization methods exist, such as the least squares stabilization first proposed for element free Galerkin (EFG) in  \cite{beissel_nodal_1996}. In this work, the square of a solution's residual is added as a stabilization term to a system's energy functional. As this term contains second order derivatives, which do not vanish at nodes, this serves to stabilize the nodal integration. Unfortunately, the computation of second derivatives, even only at nodal locations, adds considerable computational expense to this method.  

Stabilized conforming nodal integration (SCNI), introduced in \cite{chen_stabilized_2001,chen_non-linear_2002}, computes gradients via the divergence of spatially averaged strain, avoiding the computation of derivatives at nodes and also eliminating spurious oscillatory modes. 
SCNI with additional stabilization \cite{hillman_accelerated_2016} has been introduced to further increase ellipticity and hence stability of nodal integration. 
SCNI requires the construction conforming representative nodal domains, which requires additional efforts, and has thus be simplified in the form of stabilized non-conforming nodal integration (SNNI)\cite{guan_semi-lagrangian_2009,guan_semi-lagrangian_2011}. While SNNI suffers from reduced convergence rates, the addition of variational consistency terms \cite{chen_arbitrary_2013} reduces integration error to below approximation errors. 

As nodal integration requires additional stabilization terms added to the variational problem, it is challenging to adapt to the schemes to more complex multi-physics problems \cite{susuki_image-based_2024}.

\section{Interpolation-based RKPM for implementation within existing finite element software frameworks }\label{sec:Int-RKPM}

Despite the significant development of nodal integration techniques, efficient integration and assembly remains a challenge in the implementation of RKPM. This work introduces a novel computational method employing interpolation to represent RKPM bases as linear combinations of Lagrange polynomial bases. Interpolation-based RKPM (Int-RKPM) reduces assembly costs when compared to classical RKPM methods, and is implemented within existing finite element software. 

\subsection{Properties of Int-RKPM functions}

The main idea of interpolation-based methods is to replace a basis function, which for some reason may be difficult to integrate or otherwise work with, with an approximation of that function that maintains the properties of interest of the original function. This approximation is constructed of a secondary `foreground' basis $\{N_{i}\}_{i=1}^{\nu}$, where $\nu$ is the number of basis functions.  Thus, an RKPM function $\Psi_{I}$ is approximated as
\begin{align}\label{eq:approxShap}
    \Psi_{I} (\bm{x}) \sim \Hat{\Psi}_I (\bm{x}) = \sum_{j=1}^{\nu}M_{Ij} N_j(\bm{x}),   
\end{align}
where the tensor $M_{Ij}$ is referred to as the extraction operator.

The properties of a basis of interpolated functions $\{\Hat{\Psi}_I\}_{I=1}^{NP}$ will depend on the properties of both the original `background' basis $\{\Psi_I\}_{1}^{NP}$ and the foreground basis $\{N_{i}\}_{i=1}^{\nu}$. In this work, $\{N_{i}\}_{i=1}^{\nu}$ is composed of Lagrange polynomials and thus satisfies the Kronecker delta property $N_i(\bm{x}_j) = \delta_{ij}$. Due to this interpolatory property, the components of the extraction operator are easily computed by evaluating the original background function at the nodal locations of the foreground basis, 
\begin{align}\label{eq:exOp}
    M_{Ij} = \Psi_{I}(\bm{x}_j). 
\end{align}
Following \cite{fromm_interpolation-based_2023}, several properties of the interpolated basis can be derived. Provided the original basis $\{\Psi_I\}_{I=1}^{NP}$ forms a partition of unity, the interpolated basis $\{\Hat{\Psi}_I\}_{I=1}^{NP}$ will as well. The function space spanned by the interpolated basis will span polynomials of degree $\kappa$, where $\kappa =$min$\{k,n\}$, $n$ being the polynomial order spanned by the background basis and $k$ being the polynomial order spanned by the foreground basis. If both the foreground and background basis functions have local support, then the interpolated basis functions will as well.

Furthermore, if both the foreground and background basis satisfy the reproducing condition 
\begin{align}
    \sum_{I=1}^{NP} \Psi_{I}(\bm{x}) f^{\alpha}(\bm{x}_{I}) = f^{\alpha} (\bm{x}), \text{ and } \sum_{i=1}^{\nu} N_i(\bm{x}) f^{\alpha}(\bm{x}_{i}) = f^{\alpha} (\bm{x}), 
\end{align}
where $f^{\alpha}$ is a polynomial of degree $\alpha$, then the interpolated basis also satisfies the condition
\begin{align}
    \sum_{I=1}^{NP} \Hat{\Psi}_{I}(\bm{x}) f^{\alpha}(\bm{x}_{I})  = \sum_{I=1}^{NP} \sum_{j=1}^{\nu}M_{Ij} N_j(\bm{x}) f^{\alpha}(\bm{x}_{I}) =& \\
      \sum_{j=1}^{\nu} N_j(\bm{x})  \left( \sum_{I=1}^{NP}  \Psi_{I}(\bm{x}_{j})f^{\alpha}(\bm{x}_{I})\right)  = \sum_{j=1}^{\nu} N_j(\bm{x})  f^{\alpha}(\bm{x}_{j}) =& f^{\alpha} (\bm{x}). 
\end{align}

\begin{figure}
	\centering
	\includegraphics[width=1\textwidth]{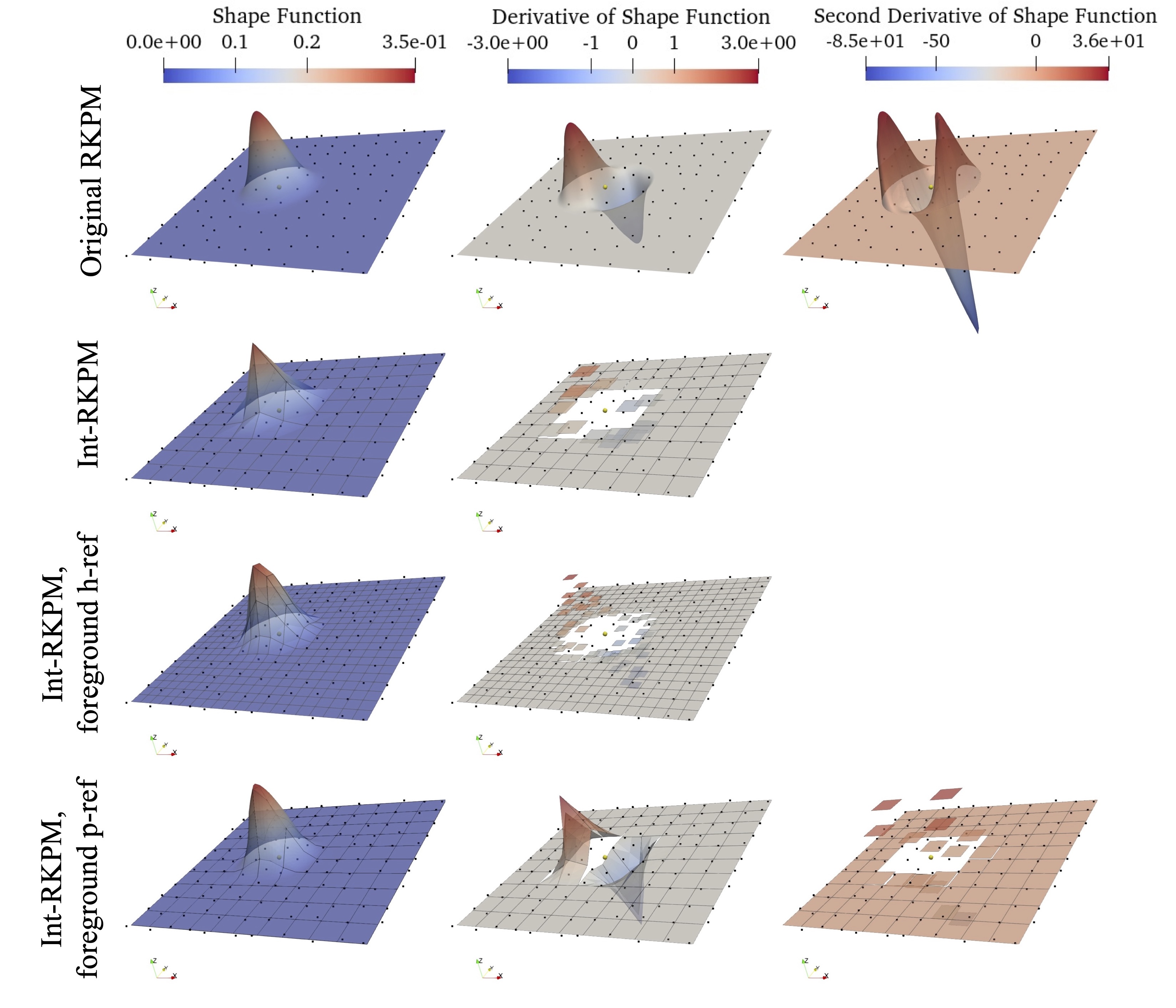}\
	\caption[ A linear RKPM shape function $\Psi_I$, using a cubic B-spline kernel is plotted along with its first derivative $\partial_{x} \Psi_I$ and second derivative $\partial_{xx} \Psi_I$.]{\label{fig:shapeFuncsLinear} A linear RKPM shape function $\Psi_I$, using a cubic B-spline kernel is plotted along with its first derivative $\partial_{x} \Psi_I$ and second derivative $\partial_{xx} \Psi_I$.  The second row shows the analogous Int-RKPM function $\Hat{\Psi}_I$ and its derivatives $\partial_{x} \Hat{\Psi}_I$ and $\partial_{xx} \Hat{\Psi}_I$, interpolated with a linear foreground function space with element size equal to the average background nodal spacing. $h$-refinement and $p$-refinement are performed on the foreground basis $N_i$, and are shown in the third and fourth rows, respectively. }
\end{figure}

\begin{figure}
	\centering
	\includegraphics[width= 1\textwidth]{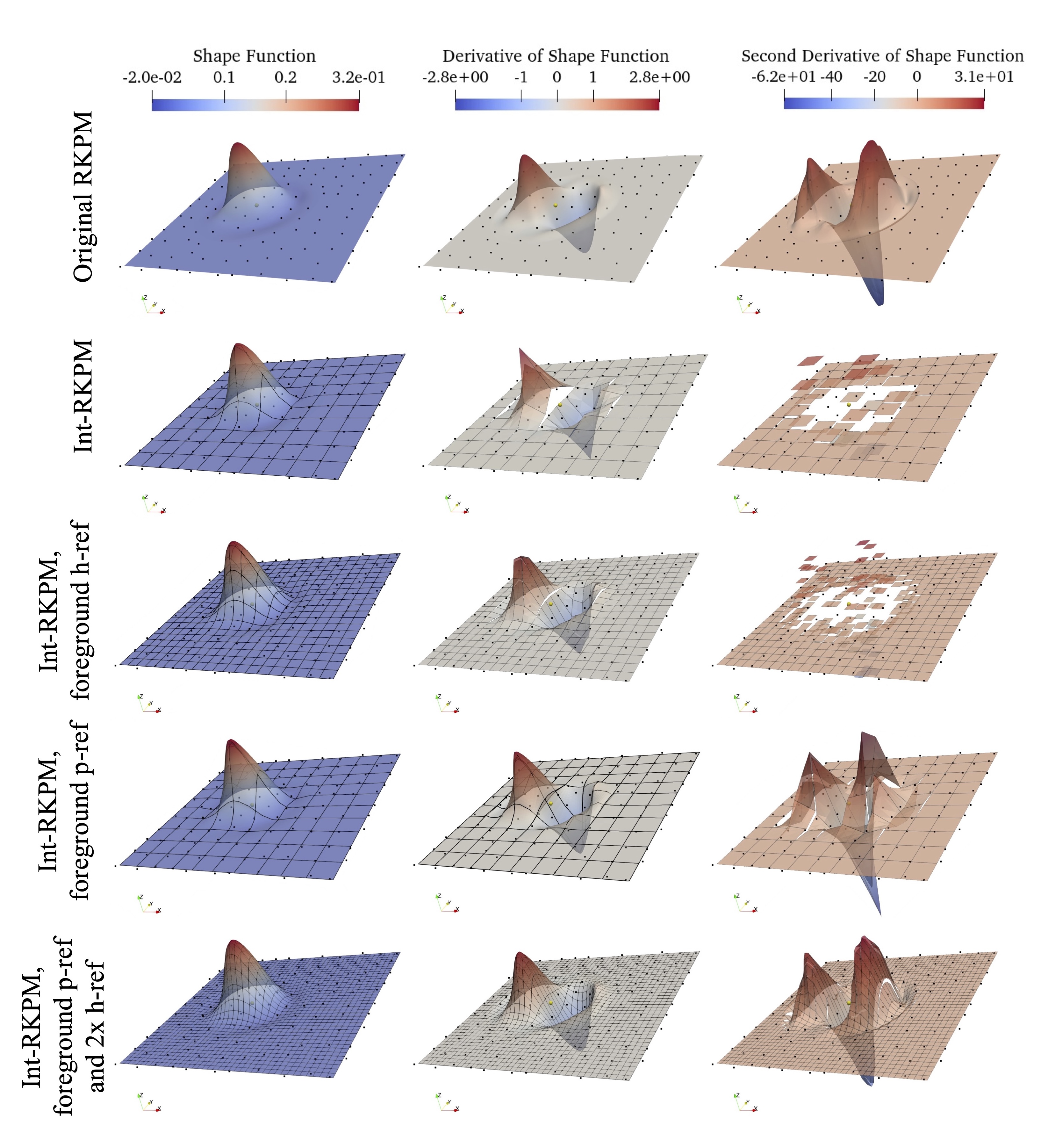}
\caption[ A quadratic RKPM shape function $\Psi_I$, using a cubic B-spline kernel is plotted along with its first derivative $\partial_{x} \Psi_I$ and second derivative $\partial_{xx} \Psi_I$. ]{\label{fig:shapeFuncsQuad}  A quadratic RKPM shape function $\Psi_I$, using a cubic B-spline kernel is plotted along with its first derivative $\partial_{x} \Psi_I$ and second derivative $\partial_{xx} \Psi_I$.  The second row shows the analogous Int-RKPM function $\Hat{\Psi}_I$ and its derivatives $\partial_{x} \Hat{\Psi}_I$ and $\partial_{xx} \Hat{\Psi}_I$, interpolated with a quadratic foreground function space with element size equal to the average background nodal spacing. $h$-refinement and $p$-refinement are performed on the foreground basis $N_i$, and are shown in the third and fourth rows, respectively. The fifth row shows the results of combined $p$-refinement and two levels of $h$-refinement, which is more accurately able to interpolate the second derivatives}
\end{figure}

The derivatives of the interpolation functions can be computed using the derivative of the foreground functions:
\begin{align}
    \partial_{x} \Hat{\Psi}_{I}(\bm{x}) = \sum_{j=1}^{\nu}M_{Ij} \partial_{x} N_j(\bm{x}). 
\end{align}
Again, provided the foreground and background bases satisfy the gradient reproducing condition
\begin{align}
    \sum_{I=1}^{NP} \partial_{x} \Phi_{I}(\bm{x}) f^{\alpha}(\bm{x}_{I}) = \partial_{x} f^{\alpha} (\bm{x}), \text{ and } \sum_{i=1}^{\nu} \partial_{x}N_i(\bm{x}) f^{\alpha}(\bm{x}_{i}) = \partial_{x}f^{\alpha} (\bm{x}), 
\end{align}
the interpolated basis will as well
\begin{align}
    \sum_{I=1}^{NP} \partial_{x}  \Hat{\Phi}_{I}(\bm{x}) f^{\alpha}(\bm{x}_{I})  = \sum_{I=1}^{NP} \sum_{j=1}^{\nu}M_{Ij} \partial_{x} N_j(\bm{x}) f^{\alpha}(\bm{x}_{I}) =& \\
      \sum_{j=1}^{\nu} \partial_{x} N_j(\bm{x})  \left( \sum_{I=1}^{NP}  \Phi_{I}(\bm{x}_{j})f^{\alpha}(\bm{x}_{I})\right)  = \sum_{j=1}^{\nu} \partial_{x}  N_j(\bm{x})  f^{\alpha}(\bm{x}_{j}) =& f^{\alpha} (\bm{x}). 
\end{align}

RKPM shape functions $\Psi$ and a series of Int-RKPM shape functions $\Hat{\Psi}$, are shown in Figures \ref{fig:shapeFuncsLinear} and \ref{fig:shapeFuncsQuad}. 

Previously published work on interpolation-based (also called extraction-based) methods, \cite{scott_isogeometric_2011,schillinger_lagrange_2016,kamensky_tigar_2019} deal primarily with exact interpolation, where the foreground basis is sufficiently refined and continuous to exactly interpolate every function in the background basis. 
In \cite{fromm_interpolation-based_2023} and \cite{fromm_interpolation-based_2024}, approximate extraction was implemented to utilize interpolation for immersed boundary analysis, where the interpolation is exact in the majority of a computational domain but approximate in the lower dimensional surface where background elements were cut. Also investigated in \cite{fromm_interpolation-based_2023} were so-called `background-unfitted' foreground meshes. In background-fitted foreground meshes, the foreground cells conform to the boundary of each background basis function support. These constraints are not present with background-unfitted foreground meshes, which nevertheless yielded optimal error convergence rates when employed to approximate PDEs. 

With RKPM and other meshfree methods the domains of support of individual basis functions are semi-random and overlapping. Thus the construction of background-fitted foreground meshes, where foreground cells would conform to the support of each RKPM function, presents a considerable challenge. As the goal of Int-RKPM is to decrease the computational cost and implementation overhead of meshfree methods, background-unfitted foreground meshes are employed. 

Int-RKPM thus represents a new iteration of background-unfitted foreground interpolation techniques. With the loosening of constraints on the foreground mesh discretization comes increased choice on the part of the user, and considerable space is dedicated in this work to quantifying the effect of foreground space on the properties of Int-RKPM. In general, the discretization used to generate the foreground basis should have average element size equal to or less than the average nodal space of the background RKPM point distribution, and should be of equal polynomial order (i.e. $k=n$). This work will discuss both foreground $h-$ refinement, where the foreground mesh is refined relative to the background RKPM discretization, and foreground $p-$ refinement, where the polynomial order of the foreground basis is increased relative to the background RKPM order of reproducibility. Both types of foreground refinement are shown in Figures \ref{fig:shapeFuncsLinear} and \ref{fig:shapeFuncsQuad}, and will be discussed in the following section.

\subsection{Using Int-RKPM bases to solve PDEs}
\label{subsec:modelProblem}
The Poisson problem is employed to model the use of Int-RKPM functions to solve PDEs. The strong form of this problem is simply:  Find $u : \Omega \rightarrow \mathbb{R}$ such that
\begin{align}
\label{eq:poisson-strong}
    \nonumber -\Delta u = f  ~~~ & \text{ in } \Omega \text{ ,} \\
            u = g  ~~~ & \text{ on } \partial \Omega \text{ ,}
\end{align}
where $g$ is Dirichlet boundary data. Nitsche's method \cite{nitsche_uber_1971} is used to enforce boundary conditions such that the approximate form can be written: Find $u^h \in \mathcal{V}^h$ such that $\forall v^h\in\mathcal{V}^h$,
\begin{align}
\int_{\Omega} \nabla u^h \cdot \nabla v^h d\Omega - \int_{\partial \Omega} \left( \nabla u^h \cdot \bm{n}\right) v^h d\Gamma \mp \int_{\partial \Omega} \left( \nabla v^h \cdot \bm{n}\right) u^h d\Gamma + \int_{\partial \Omega} \frac{C_\text{pen}}{h} u^h v^h d\Gamma \nonumber \\
= \int_{\Omega} f v^h d\Omega \mp \int_{\partial \Omega} \left( \nabla v^h \cdot \bm{n}\right) g d\Gamma + \int_{\partial \Omega} \frac{C_\text{pen}}{h} g v^h d\Gamma\text{ ,}\label{eq:poisson-disc-generic}
\end{align}
where $\bm{n}$ is the normal vector, $C_\text{pen}$ is a user defined constant, and $h$ is understood to be the background function nodal spacing. This work employs the symmetric Nitsche's method, such that the $\mp$ signs in equation \ref{eq:poisson-disc-generic} are taken as $-$. The left and right hand side of equation \ref{eq:poisson-disc-generic} can be grouped into a bilinear and linear form, such that, 
\begin{align}
    a (u^h, v^h) = L(v^h). \label{eq:bilinearform}
\end{align}

For classic RKPM, the function space $\mathcal{V}^h = \mathcal{V}_{RKPM} = \text{span} \{\Psi_I\}$, where each shape function $\Psi_I$ can be computed with Equation \ref{eq:RKPMShapeCOmp} . Functions $u^h_{RKPM} \in \mathcal{V}_{RKPM}$ are then expressed as
\begin{align}
    u^h_{RKPM} (\bm{x}) = \sum_{I} \Psi_I(\bm{x}) u_I. 
\end{align}
These shape functions can be integrated either using Gauss-quadrature schemes or some flavor of stabilized nodal integration, to assemble the linear system 
\begin{align}
\bm{K} \bm{d} = \bm{F} \label{eq:matrix}
\end{align}
where $K_{IJ}=  a (\Psi_I, \Psi_J) $, $d_I = u_I$, and $F_I = L(\Psi_I)$ are defined in terms of the bilinear and linear forms. 

With interpolation-based RKPM the function space $\mathcal{V}^h = \mathcal{V}_{Int-RKPM} = \text{span} \{\Hat{\Psi}_I\}$, where each shape function $\Hat{\Psi}_I$ is an approximation of the corresponding function $\Psi_I$, as in equation \ref{eq:approxShap}. With the extraction operator given in equation \ref{eq:exOp}, functions $u^h_{Int-RKPM} \in \mathcal{V}_{Int-RKPM}$ are expressed as
\begin{align}
    u^h_{Int-RKPM} (\bm{x})= \sum_{I} \Hat{\Psi}_I(\bm{x}) \Hat{u}_I = \sum_{I} \sum_{j} M_{Ij} N_j(\bm{x}) \Hat{u}_I. 
\end{align}

To assemble the linear system in equation \ref{eq:matrix}, 
allow $\Hat{K}_{IJ}=  a (\Hat{\Psi}_I,\Hat{\Psi}_J)$, $\Hat{d}_I = \Hat{u}_I$, and $\Hat{F}_I = L(\Hat{\Psi}_I) $. The stiffness matrix and force vector are thus 
\begin{align} \label{eq:IntStiffness}
    \Hat{K}_{IJ}=  a (\Hat{\Psi}_I,\Hat{\Psi}_J) = M_{Ik} a (N_k,N_l)  M_{lJ }\text{ and }\Hat{F}_I = L(\Hat{\Psi}_I) = M_{Ik} L(N_k).
\end{align}
The quantities $a(N_k,N_l) $ and $L(N_k)$ can be easily computed with the foreground Lagrange polynomial basis $\{N_{i}\}_{i=1}^{\nu}$. The only evaluation of the original RKPM basis resides within the computation of the extraction operator, see equation \ref{eq:exOp}. This feature of Int-RKPM allows existing finite element software to be augmented to use meshfree methods with the introduction of an extraction operator.

\subsubsection{Implementation through open-source finite element software}

Efficient integration techniques have been developed for meshfree methods, but these techniques are difficult to implement within existing finite element codes. Open source implementations of RKPM \cite{huang_rkpm2d_2020} and other meshfree methods \cite{ramachandran_pysph_2021} rely on considerable custom software packages.

As Int-RKPM replaces the integration of RKPM basis functions with integration of classic Lagrange polynomials, it can be naturally implemented within existing FEM software. Following previous works on interpolation-based IGA \cite{kamensky_tigar_2019, kamensky_open-source_2021} and immersogeometric methods \cite{fromm_interpolation-based_2023, fromm_interpolation-based_2024}, this work employs the popular open source software FEniCSx \cite{baratta_dolfinx_2023} to illustrate the potential of Int-RKPM.

\subsection{Numerical comparison of int-RKPM with classic RKPM}
\label{subsec:Poisson}

To demonstrate the efficacy of Int-RKPM a simple Poisson problem on a unit square was modeled. A suite of RKPM point sets was used to investigate error convergence rates. Each point set was created by perturbing a uniform grid of points in a jittered grid \cite{cook_stochastic_1986} with initial uniform spacing of $h$. To construct the jitter grid, each coordinate was perturbed such that $(x,y)\rightarrow (x,y) + \epsilon h (\eta_x, \eta_y) $. A number $\eta \in(-1\,, 1)$ was randomly generated for each point and coordinate and $\epsilon$ is referred to as the perturbation parameter. The accuracy of integration schemes employed in classic RKPM is shown to be linked to to the uniformity of points \cite{chen_arbitrary_2013} thus multiple perturbation parameters $\epsilon \in [0,1]$ were investigated. However, provided points did not overlap, the performance of Int-RKPM was not dependent on the point distribution. All RKPM point sets in work thus use $\epsilon = 0.5$. 

RKPM bases were defined on these point sets using cubic B-spline kernels and circular supports. The normalized support size of each function was $n+1$, where $n$ is the polynomial order of reproducibility. Both linear ($n=1$) and quadratic ($n=2$) RKPM bases were investigated. 

Interpolation was performed using continuous-Galerkin type Lagrange polynomial foreground bases defined upon a suite of uniform quadrilateral boundary conforming meshes. The effects of both $p-$ and $h-$refinement on the foreground basis were investigated. Initial foreground element sizes are proportional to background nodal spacing. With $h-$refinement, the ratio of foreground element size to average RKPM nodal spacing was decreased to 1:2 (for 1x $h-$ref) and 1:4 (for 2x $h-$ref). The polynomial order of the foreground basis, $k$, is set initially to $n$, the polynomial order of reproducibility of the background basis. With $p-$refinement, $k = n + p$ where $p$ is the level of refinement.

Int-RKPM is compared with a classic implementation of RKPM where integrals are computed using a Gauss-Quadrature grid. The grid size is equal to the average nodal spacing of the RKPM pointset. A $6\times6$ Gauss scheme is used for the linear RKPM basis, and an $8\times8$ Gauss scheme is used for the quadratic RKPM basis.

The Poisson problem as introduced in \ref{subsec:modelProblem} with Equation \ref{eq:poisson-strong} is modeled. Using the method of manufactured solutions the source term and boundary data are set to $f = - \Delta u_{ex}$ and $g= u_{ex}$, respectively,  where 
\begin{align}
    u_{ex} = \sin(0.1x + 0.1)\sin(0.1y + 0.1). 
\end{align}

\begin{figure}
	\centering
	\includegraphics[width= \textwidth]{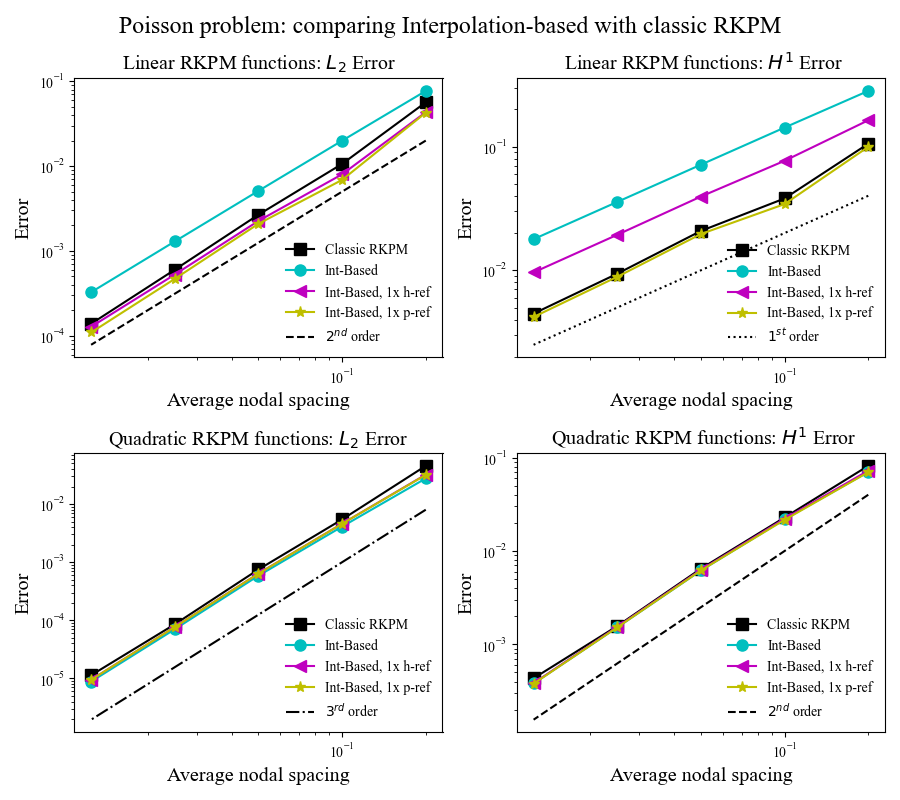}\
	\caption{\label{fig:poisson}Convergence data for the Poisson problem, comparing Int-RKPM with classic RKPM.}
\end{figure}

For the linear basis, results shown in the first row of Figure \ref{fig:poisson}, interpolation-error increases error magnitude but does not impact the convergence rate of Int-RKPM. Foreground $h-$refinement reduces error magnitude to below that of classic RKPM for the $L_2$ error norm, and reduces the magnitude marginally for the $H_1$ error norm. A level of foreground $p-$refinement is sufficient to match the error magnitude of Int-RKPM to classic RKPM for this example. 

This behaviour is qualitatively explained by comparing the derivatives of Int-RKPM in Figure \ref{fig:shapeFuncsLinear}. The second row of Figure \ref{fig:shapeFuncsLinear} corresponds to the blue (circle) line in Figure \ref{fig:poisson}, where the first derivative is interpolated as the linear combination of relatively coarse piecewise constant foreground functions. When foreground $h-$ refinement is applied for the third row of Figure \ref{fig:shapeFuncsLinear}, the approximation is improved, resulting in the decreased error magnitude shown by the pink (triangle) line in the plot in Figure \ref{fig:poisson}. However, it is not until foreground $p-$refinement is applied to in the fourth row of \ref{fig:shapeFuncsLinear}, and the interpolation of the derivative is given by pieceswise linear functions, that the interpolation error is decreased sufficiently to reproduce the classic RKPM results. This is shown by the overlap of the yellow (star) line and black (square) line in the plot in Figure \ref{fig:poisson}. However, it must be emphasized that no amount or type of local foreground refinement is required for the Int-RKPM method to exhibit optimal convergence rates.

For the quadratic basis, the second row of plots in Figure \ref{fig:poisson}, the results produced by Int-RKPM are practically identical to classic RKPM. Additional foreground refinement neither increases the rate of convergence nor decreases the magnitude of the error. These results can be qualitatively understood by analyzing the quadratic Int-RKPM derivatives in Figure \ref{fig:shapeFuncsQuad}. The second row of the figure shows the Int-RKPM functions without any foreground refinement. The first derivative of the function is interpolated with piecewise linear functions. The improvements to the interpolation of the derivative are minimal with the application of foreground $h-$refinement, seen in the third row, and foreground $p-$refinement, seen in the fourth row. It is worth noting that in all cases involving foreground h- and p-refinements, the total number of degrees of freedom for the background RKPM basis remains unchanged. Consequently, the dimensions of the stiffness matrix and force vector do not increase.

\subsection{Reducing computational costs with high-order derivatives }
\label{subsec:biharmonic}

Numerous problems in engineering applications are described by PDEs containing higher order derivatives. Notably, the equations characterizing Kirchhoff-Love shells contain fourth order derivatives and must be approximated by functions belonging to an $H^2$-conforming basis. RKPM bases can easily satisfy this property, but due to the relative complexity in computing high order basis function derivatives, particle methods have not been widely applied to Kirchhoff-Love shell problems.

The computation of derivatives of an interpolated RKPM basis can be done with standardized finite element operations. To demonstrate the efficacy of interpolated-RKPM with non-conforming function spaces, the biharmonic problem is tested, with strong form: 
 Find $u : \Omega \rightarrow \mathbb{R}$ such that 
\begin{equation}
    \Delta^2 u = f\text{ ,}
\end{equation}
with boundary conditions
\begin{align}
    u = g &\text{ on } \partial\Omega\text{ ,} \\
    \grad u \cdot \boldsymbol{n} = h  &\text{ on } \partial\Omega\text{ ,}
\end{align}
where $f:\Omega\to\mathbb{R}$ is a given source term and $g:\partial \Omega\to\mathbb{R}$ and  $h:\partial \Omega\to\mathbb{R}$ are boundary data.
Employing the method of manufactured solutions, $f = \Delta^2 u_{ex}$, $g = u_{ex}$, and $h = \grad u_{ex} \cdot \boldsymbol{n}$. For this numerical example the same exact solution is employed, 
\begin{equation}
    u_{ex}(x,y) = \sin (0.1y + 0.1)  \sin (0.1x +0.2) .
\end{equation}
Boundary conditions are weakly enforced using Nitsche's-like residual terms, such that the weak form can be written as: Find $u^h\in\mathcal{V}^h$ such that, $\forall v^h\in\mathcal{V}^h$, 
\begin{align} 
    \nonumber \int_{\Omega}\Delta u^h\Delta v^h\,d\Omega + \int_{\partial\Omega}\grad \Delta u^h \cdot \boldsymbol{n}v^h - \Delta u^h \grad v^h \cdot  \boldsymbol{n}\,d\Gamma & \\
    \nonumber  + \int_{\partial\Omega}(\grad \Delta v^h) \cdot \boldsymbol{n}(u^h - u_{ex})  - \Delta v(\grad u^h \cdot  \boldsymbol{n} - \grad u_{ex}  \cdot  \boldsymbol{n})\,d\Gamma  & \\
    + \int_{\partial\Omega}\dfrac{\alpha}{h^3} (u^h - u_{ex})v^h + \dfrac{\beta}{h} (\grad u^h \cdot  \boldsymbol{n} - \grad u_{ex}\cdot  \boldsymbol{n})\grad v^h\cdot  \boldsymbol{n}\,d\Gamma  & = \int_\Omega f v^h\,d\Omega\text{ ,}\label{eq:biharmonic-nitsche}
\end{align}
where $\alpha > 0$ and $\beta > 0$ are user specified constants. For the computations of this paper, $\alpha=\beta=10$. 

\begin{figure}
	\centering
	\includegraphics[width= \textwidth]{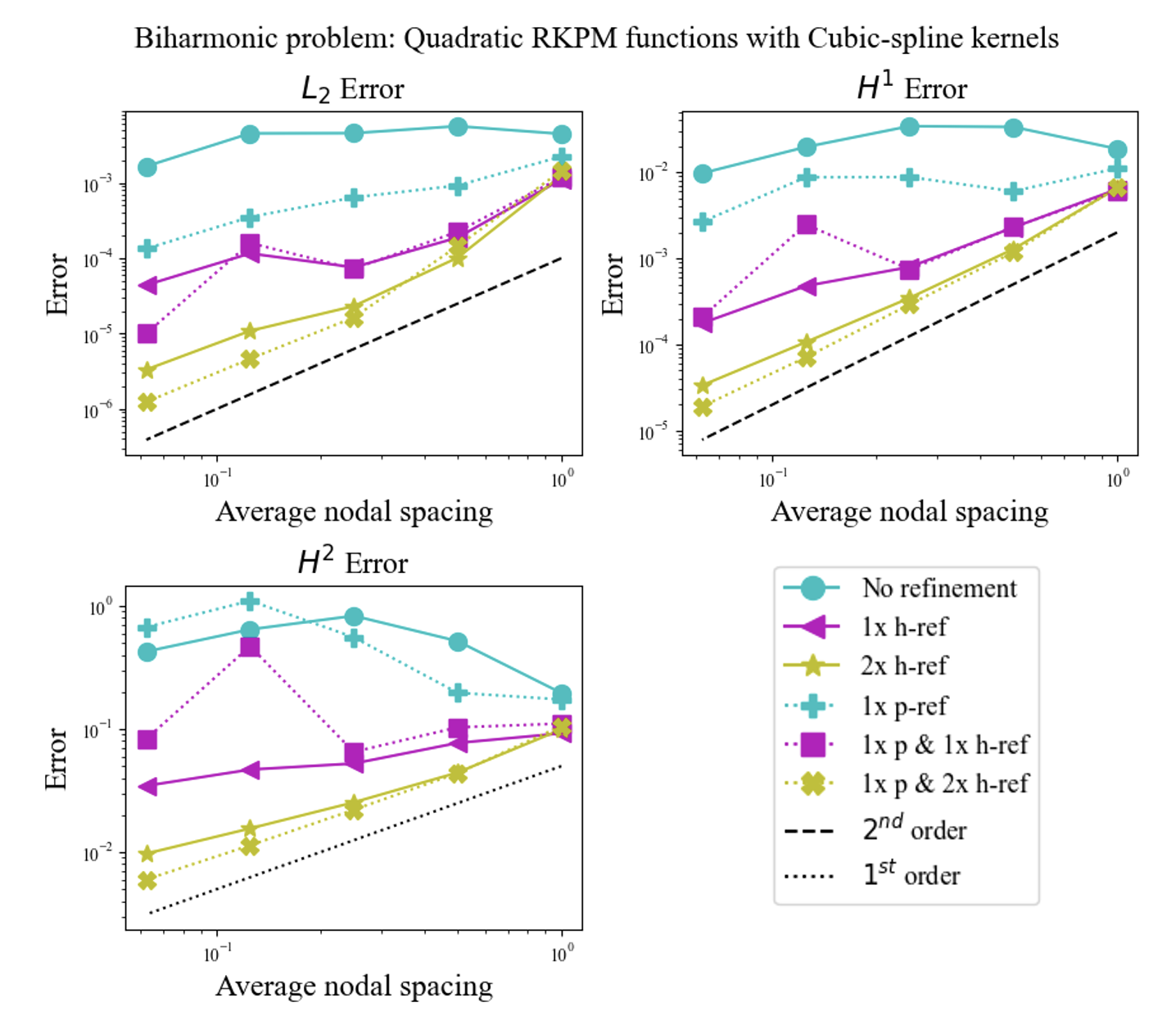}\caption{\label{fig:biharmonic}Convergence data for the biharmonic problem approximated with Int-RKPM, comparing various levels of foreground $h-$ and $p-$ refinement.}
\end{figure}

The same suite of background RKPM basis functions from the previous example are used along with uniform quadrilateral foreground meshes. Numerical tests investigated $L_2$, $H^1_0$, and $H^2_0$ error convergence rates with concurrent foreground and background refinement. Results are shown in Figure \ref{fig:biharmonic}.

As in the lower order derivative Poisson problem investigated in subsection \ref{subsec:Poisson}, some additional foreground refinement is required to reduce interpolation error. This is because with this example, the method must not only interpolate the first derivatives of the original RKPM function, but also the second derivatives. As seen in the first rows of Figures \ref{fig:shapeFuncsLinear} and \ref{fig:shapeFuncsQuad}, the second derivatives of RKPM functions are difficult to fully resolve. This, in addition to the high computational cost of computing higher order derivatives, is a reason why meshfree methods are under utilized for high order PDEs like the biharmonic problem. 

While the results of the previous section showed that $p-$refinement of the foreground mesh improved the approxiation potential of Int-RKPM, $p-$refinement alone is not sufficient to achieve optimal error convergence rates for this PDE, as seen by the solid colored pink blue and green lines. And unlike the Poisson problem, the original Int-RKPM without foreground refinement (the solid pink line) not only has a larger error magnitude, but the errors do not converge at all. In this case, foreground $h-$refinement is required to achieve optimal rates. With two levels of $h-$refinement (such that the ratio of foreground element size to background nodal spacing is 1:4 ) optimal rates are observed. This can be qualitatively understood from the third column of images in Figure \ref{fig:shapeFuncsQuad}. While $p-$refinement alone was sufficient to reduce the error magnitude for the linear Poisson problem, the second derivative of the interpolated basis function still contains major discontinuities that required both $p-$ and $h-$refinement. 

\section{Achieving high-order accuracy with local foreground  refinement}\label{sec:localRef}

Int-RKPM can also be applied to problems of more complicated geometric domains than those illustrated in the previous two sections. RKPM, and other meshfree methods, do not require high quality boundary conforming meshes to approximated curved or otherwise complex geometries, but in general require a boundary conforming geometric discretization to perform integration. This is likewise the case for Int-RKPM where the integration discretization is in fact a mesh used to create a foreground function space. The following section considers a square domain with a quarter-circle hole.

\subsection{Quadtree local refinement strategy}
\begin{figure}[t]
	\centering
	\includegraphics[width= \textwidth]{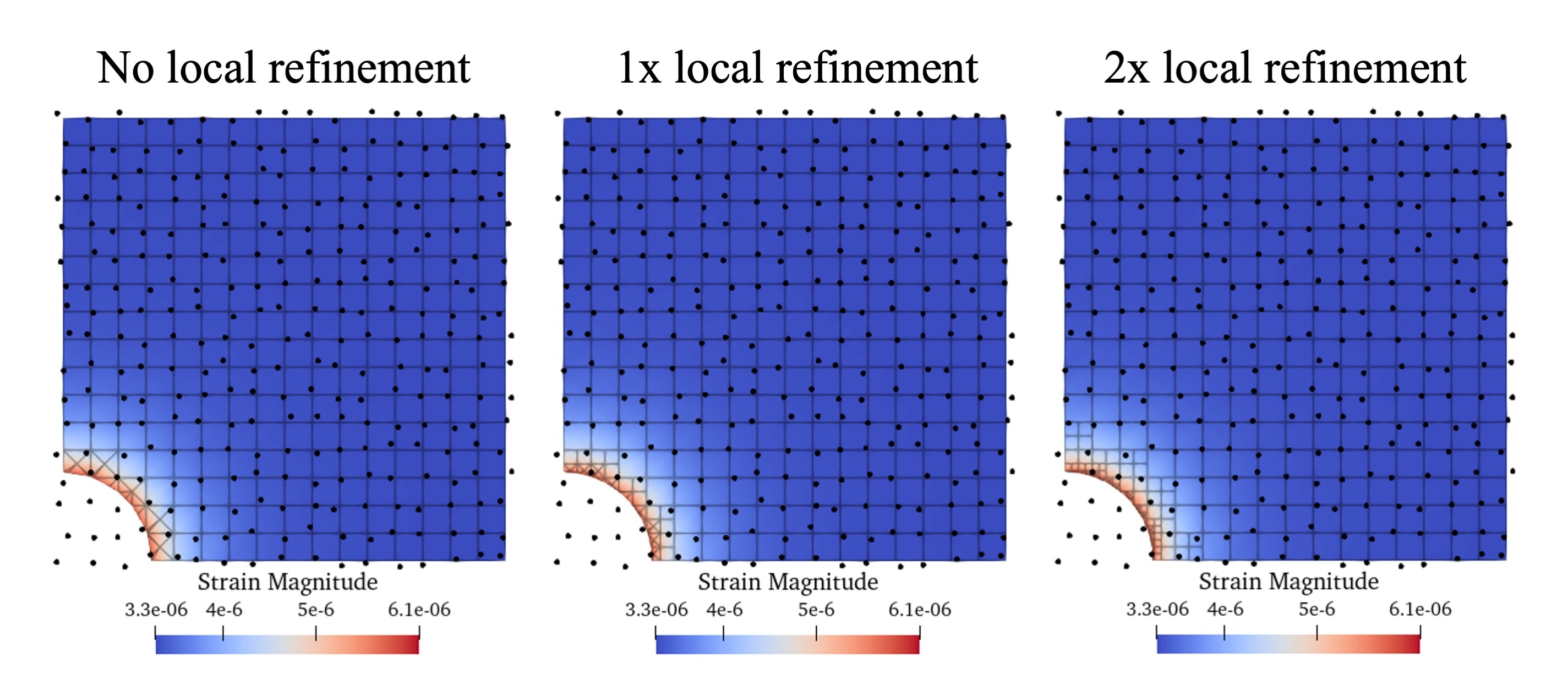}\
	\caption{\label{fig:HIPPic} Local foreground refinement is performed to increase geometric resolution for problems with approximate interface representation. The foreground meshes are shown along with the RKPM point discretization, which is the same for each level of local refinement.}
\end{figure}

The boundary conforming foreground meshes used for this domain and shown in Figure \ref{fig:HIPPic} were generated using the open source software MORIS \cite{maute_moris_2023}. MORIS is an immersed finite element code that utilizes level-set geometry descriptions to enrich and integrate both classic Lagrange and isogeometric basis functions. It also provides functionality to rapidly generate boundary conforming meshes with its level-set geometry processing sublibrary. In 2D, the algorithm begins with a uniform tensor-product grid, and then triangulates the quadrilateral cells intersected by the level-set isocountour indicating the domain boundary. After triangulation a root finding algorithm is used to locate the intersections of element edges and the isocounter. Using the intersection points as new nodes, the cells are further triangulated. The analogous 3D process begins with tetrahedron cells that are decomposed into a boundary conforming mixed hexahedron-tetrahedron mesh. 

Local foreground refinement is performed to improve geometric resolution and reduce geometry error. The algorithm to implement local refinement employs quadtree or (or octree in 3D) refinement of cells in the intial background grid prior to their triangulation (or subdivision into hexahedrons). If additional levels of refinement are required, the cells of the newly refined grid intersected by the level-set isocontour can be further identified and subdivided. For this example, three sets of foreground meshes are utilized, with no, one level, and two levels of local refinement, examples of which are shown in Figure \ref{fig:HIPPic}. 

The meshes generated for this purpose are not suitable for classical FEM as they are poorly conditioned and may contain elements with arbitarily large aspect ratios. Additionally, the meshes generated with local refinement contain hanging nodes which are not supported by most FEM software. While the poor conditioning of the elements on the foreground mesh presents no challenge for interpolation based methods, the hanging nodes are a different story. Function spaces defined on meshes with hanging nodes are $C^{-1}$ continuous, and are not suitable for approximating PDEs with standard Galerkin's method. While the previous numerical example in Section \ref{subsec:biharmonic} utilized a non-conforming function space, it was shown that foreground refinement was required to reduce interpolation error and achieve optimal convergence rates. The problem presented by a  $C^{-1}$  foreground mesh is similar, but cannot be overcome with foreground refinement alone. 

Instead, a novel `double' interpolation strategy is employed, refered to as D-Int-RKPM. A new $C^{0}$ `midground' mesh and function space $\mathcal{V}^{mg} = \text{span} \{ N^{mg}_{i}\}_{i=1}^{\mu}$ is introduced. The midground mesh is selected such that the foreground mesh is `midground fitted', meaning the domain of each cell on the midground can be exactly represented as the union of cells on the foreground mesh. This minimizes interpolation error between the midground and foreground mesh and ensures the continuity of the interpolated basis. In this case, the midground mesh is taken as the initial grid used in the foreground mesh generation process. 

The RKPM basis is first interpolated with functions from the midground space, and then interpolated again onto the foreground space
\begin{align}
    \Hat{\Psi}_I(\bm{x}) = \sum_{j}^{\mu} M^1_{Ij} N^{mg}_{j}(\bm{x}) = \sum_{j}^{\mu} \sum_{k}^{\nu} M^1_{Ij} M^2_{jk} N_{k} (\bm{x}),
\end{align}
where 
\begin{align}
    M^1_{Ij}  = \Psi_{I}(\bm{x}^{mg}_{j}) \qquad \text{and} \qquad 
    M^2_{ik}  = N^{mg}_{i}(\bm{x}_{j})
\end{align}
are the new extraction operators, with $\bm{x}^{mg}$ the nodal coordinates of the new midground basis. 

\begin{figure*}
    \centering
	 \begin{subfigure}[b]{\linewidth}
	 \centering
    \includegraphics[width=\linewidth]{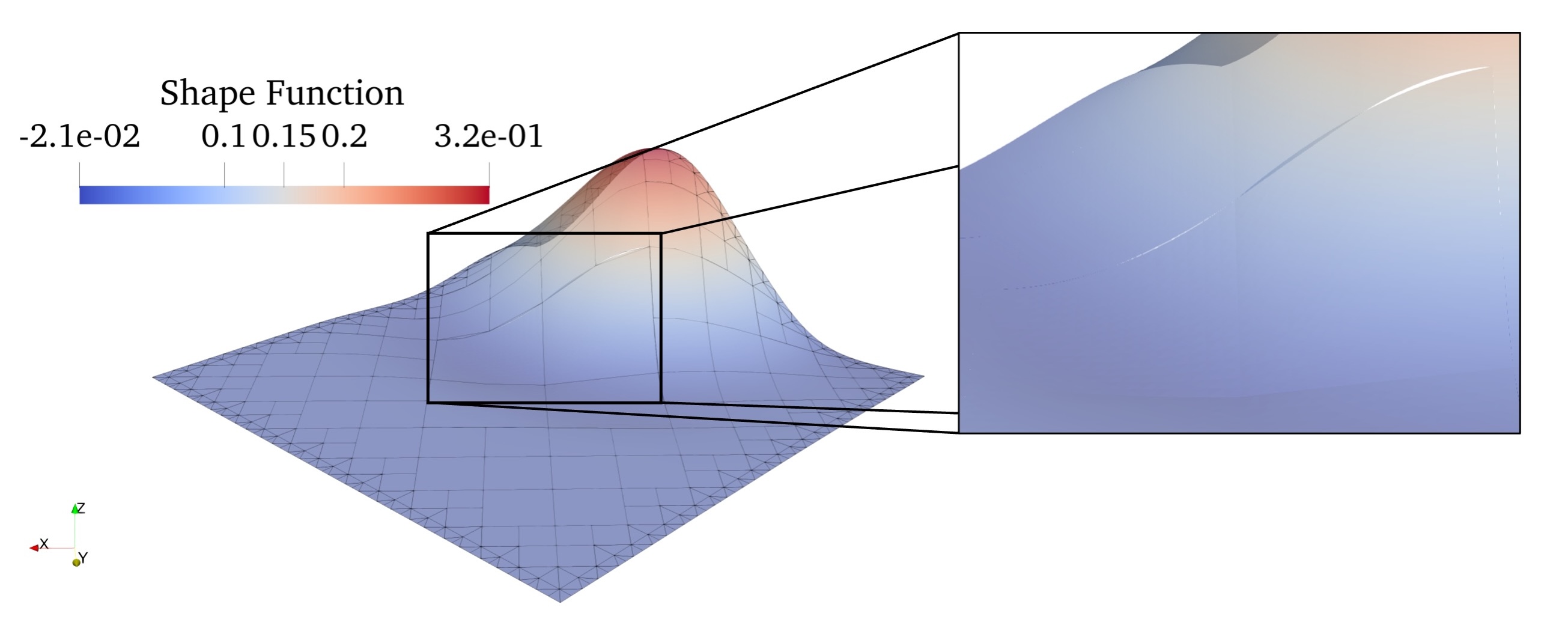} 
    \caption{A quadratic shape function is interpolated with a foreground function space defined on a mesh with hanging nodes, resulting in discontinuities.}
  \end{subfigure}
	 \begin{subfigure}[b]{1\linewidth}
	 \centering
    \includegraphics[width=\linewidth]{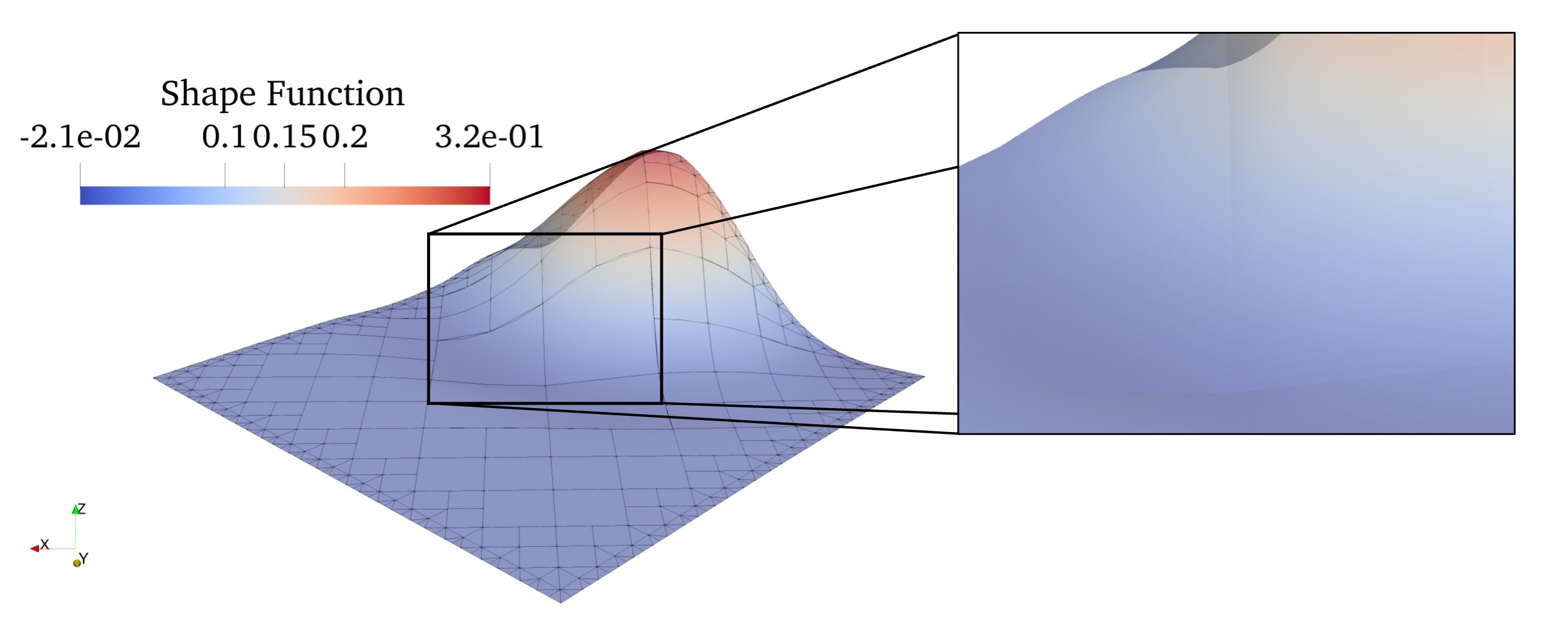} 
    \caption{The same shape function is first interpolated with a $C^0$ midground function space, and then with the $C^{-1}$ foreground space, eliminating discontinuities}
  \end{subfigure}

    \caption{Visualizations of an RKPM shape function interpolated with a $C^{-1}$ function space, and in the two step process referred to as double interpolation (D-Int) using an intermediary $C^{0}$ midground space. }
    \label{fig:hangingNodes}
\end{figure*}

The impact of this double interpolation strategy on an individual shape function is illustrated in Figure \ref{fig:hangingNodes}.

\subsection{Linear elasticity modeled on a curved geometry}\label{subsec:holeInPlate}

In this example linear elasticity is modeled in a plate with a hole. Navier's equations for linear elasticity are as follows: Find $\bm{u}: \Omega \rightarrow  \mathbb{R}^2 $
\begin{align}
   \grad \cdot \bm{\sigma}(\bm{u}) = 0
\end{align}
with Cauchy stress tensor 
\begin{align}
    \sigma_{ij} (\bm{u}) = 2 \mu \varepsilon_{ij} (\bm{u}) + \lambda \varepsilon_{kk}(\bm{u})\bm\delta_{ij},
\end{align}
where the $\varepsilon(\bm{u})_{ij} =\dfrac{1}{2}( u_{i,j} + u_{j,i})$ is the strain and $\mu=151$GPa and $\lambda=65.9$GPa are the Lam\'{e} parameters. Traction and symmetry boundary conditions are  imposed in the plate with a hole problem as demonstrated in Figure \ref{fig:HIPPic}:
\begin{align}
    \bm{u}\cdot \textbf{n} = 0\quad&\text{on}~\Gamma_\text{sym}\text{ ,}\\
    (\bm{I}-\boldsymbol{n}\otimes\bm{n})\cdot (\bm{\sigma}\cdot\textbf{n}) = \bm{0}\quad&\text{on}~\Gamma_\text{sym}\text{ ,}\\
    \bm{\sigma}\cdot\textbf{n}= \bm{t}\quad&\text{on}~\Gamma_t\text{.}
\end{align}

Using Nitsche's method to enforce symmetry conditions, Navier's equation is discretized as: Find $\bm{u}^h\in\bm{\mathcal{V}}^h = [\mathcal{V}^h, \mathcal{V}^h]$ such that, $\forall\bm{v}^h\in\bm{\mathcal{V}}^h$,
\begin{align}
    \nonumber &\int_\Omega \bm{\sigma}(\bm{u}^h):\nabla\bm{v}^h\,d\Omega\\
    \nonumber &\quad - \int_{\Gamma_\text{sym}}(\bm{u}^h\cdot\textbf{n})\textbf{n}\cdot\bm{\sigma}(\bm{v}^h)\cdot\textbf{n}\,d\Gamma - \int_{\Gamma_\text{sym}} (\bm{v}^h\cdot\textbf{n})\textbf{n}\cdot\bm{\sigma}(\bm{u}^h)\cdot\textbf{n}\,d\Gamma \\
    &\quad + \int_{\Gamma_\text{sym}}\frac{\beta \mu}{h}(\bm{u}^h\cdot\textbf{n})\textbf{n}\cdot\bm{v}^h\,d\Gamma = \int_{\Gamma_t}\bm{t}\cdot\bm{v}^h\,d\Gamma\text{,}
    \end{align}
where $\beta > 0$ is the penalty parameter associated with Nitsche's method, is set to $10$.
Exploiting symmetry, a $5\times5$ square with a quarter circle of radius $R=1$ was modeled. For an infinite plate under an applied equal-biaxial strain, the resulting stress tensor components can be computed using Kirsch's equations as  
\begin{align}
    \sigma_{rr}  = \sigma_{\infty} \left( 1 - \left(\dfrac{r}{R}\right)^2 \right), \qquad&\text{and}\qquad \sigma_{\theta\theta}  = \sigma_{\infty} \left( 1 + \left(\dfrac{r}{R}\right)^2 \right). 
\end{align}
The exact solution for the stress field, converted to Cartesian coordinates is applied as traction $\bm{t} = \bm{\sigma}_{ex}\cdot \textbf{n}$. 

Jittered grid RKPM pointsets with $\epsilon = 0.5$ were generated without regard to the geometry of the hole/inclusion. A series of pointsets with average nodal spacing  $h$  were used to generate convergence data, with $h\in 0.625\{ 1, 0.5, 0.25, 0.125, 0.06125 \}$. The exact stress field, domain geometry, and the coarsest point set are all shown in Figure \ref{fig:HIPPic}. This Figure also illustrates the effects of both local refinement and double interpolation on the approximation PDEs is shown in Figure \ref{fig:HIPPic}. 

\begin{figure}
	\centering
	\includegraphics[width= \textwidth]{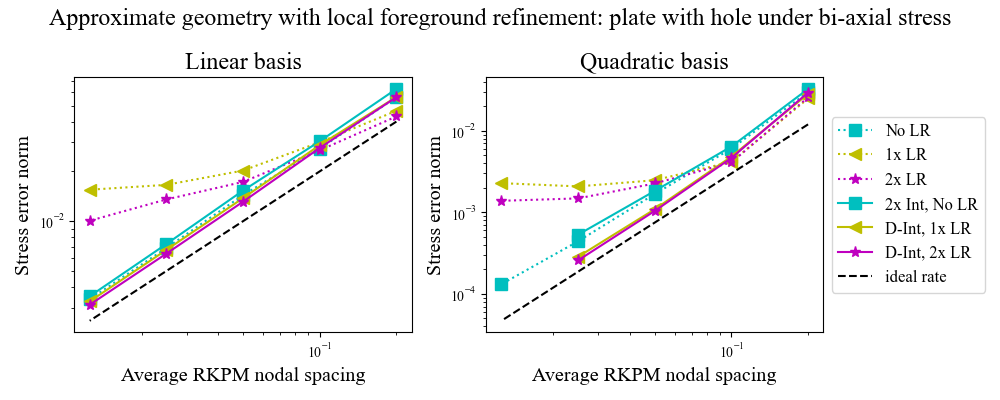}\
	\caption{\label{fig:HIPPlot}Convergence data for Int-RKPM applied to the linear elasticity problem of a plate with hole. Normal interpolation is used for the dotted lines, and double interpolation used for the solid lines}
\end{figure}

Figure \ref{fig:HIPPlot} gives the results of this PDE approximation and shows the effects of the double interpolation technique on the stress error norm. The dotted lines depict results from the original interpolation strategy, while solid lines shown results with double interpolation. With the original interpolation strategy, the discontinuities in the locally refined interpolated bases degrade error convergence rates. With the introduction of doubled interpolation, local refinement is sufficient to reduce geometric error and achieve ideal error convergence rates for both linear and quadratic interpolated RKPM functions.

\section{Heaviside enrichment for modeling multi-material domains}\label{sec:heaviside}

\subsection{Three-material square model problem}\label{subsec:3mat}

Heaviside enrichment is employed to modify RKPM for multi-material problems, following the enrichment scheme employed for immersed IGA in \cite{fromm_interpolation-based_2024}. Material subdomains are described with characteristic functions $\psi^m$, 
\begin{equation} \label{eq:charFunc}
    \psi^m(\bm{x}) = \begin{cases}
        1, \text{ if } \bm{x} \in \Omega^m\\
        0, \text{ else.}
    \end{cases}
\end{equation}
Each basis function is inspected to determine if its support is intersected by a material interface. If so, it is enriched. The enriched basis functions are expressed as 
\begin{align}\label{eq:enrichB}
    \Psi^m_I(\bm{x}) = \psi^m (\bm{x})\Psi(\bm{x}), \ \ \forall \ m \in \{1, ..., L\},
\end{align}
where $L$ is the number of material subdomains. The original basis function can be recovered by the summation of these enriched functions,
\begin{align}\label{eq:enrichBSum}
    \Psi_I(\bm{x}) = \sum_{m} \Psi^m_I(\bm{x})
\end{align}
This section implements global enrichment as opposed to function-wise enrichment \cite{schmidt_extended_2023}, requiring sufficiently small supports of basis functions in the vicinity of the interface.  

\begin{figure}
    \centering
	 \begin{subfigure}[b]{0.3\linewidth}
	 \centering
    \includegraphics[width=\linewidth]{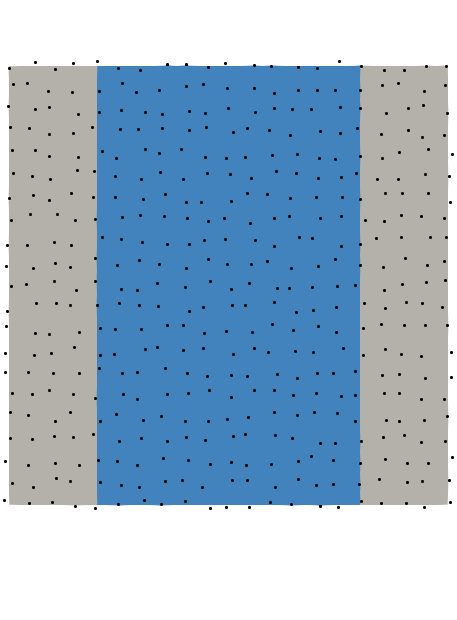}
    \caption{}
  \end{subfigure}
	 \begin{subfigure}[b]{0.3\linewidth}
	 \centering    \includegraphics[width=\linewidth]{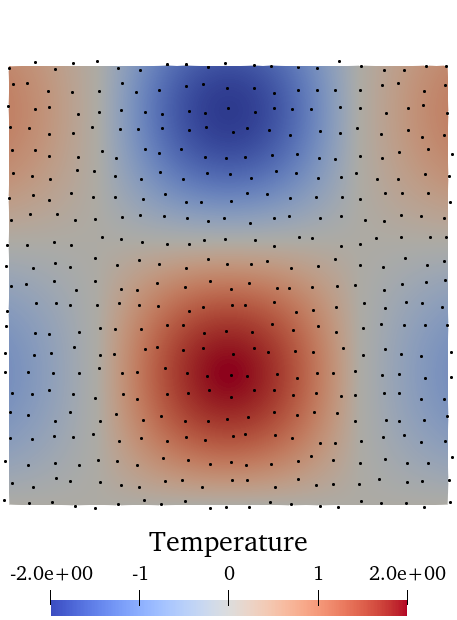}
    \caption{}
  \end{subfigure}
  \begin{subfigure}[b]{0.3\linewidth}
	 \centering
    \includegraphics[width=\linewidth]{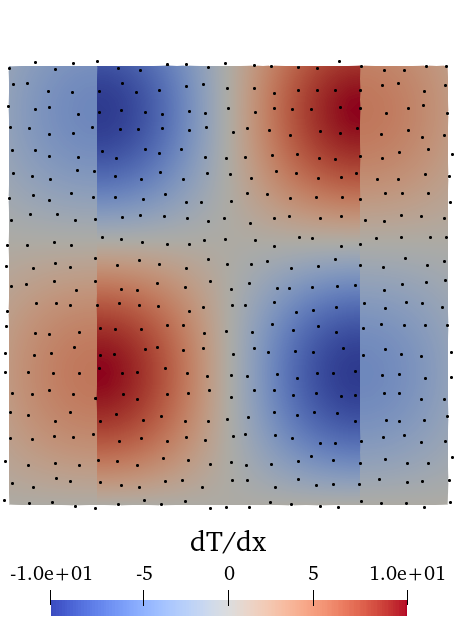}
    \caption{}
  \end{subfigure}
    \caption{The subdomains, exact solution $T_{ex}$, and the x-derivative of the exact solution. Black dots show the nodal locations of the RKPM basis functions.}
    \label{fig:threeMat}
\end{figure}

For enriched Int-RKPM, the material characteristic functions are discretized with a piece-wise constant function space defined on the foreground integration mesh. This enrichment scheme is tested on a multi-material domain using Poisson's problem. A 1x1 unit square is divided into 3 subdomains with interfaces at $x=0.2$ and $x = 0.8$, such that the subdomains are defined
\begin{align} \label{eq:analyticDomain}
    \nonumber\bm{x} \in \Omega^1, \,  0.0 \leq x \leq 0.2 \\
    \bm{x} \in \Omega^2, \, 0.2 < x \leq 0.8 \\
    \nonumber \bm{x} \in \Omega^3, \, 0.8 < x \leq 1.0 .
\end{align}
The subdomains are assigned material properties $\kappa(\bm{x}) = \kappa^m$, $\bm{x} \in \Omega^m$, with $\kappa^1 = 1.0$, $\kappa^2 = 0.5$, and $\kappa^3 = 1.0$, and are shown in Figure \ref{fig:threeMat}(a).

A source term $f:\Omega\rightarrow\mathbb{R}$ and Dirichlet boundary data $\overline{T}:\partial \overline{\Omega}\rightarrow\mathbb{R}$ on $\Gamma_{\overline{T}} \! \subset \!\partial \overline{\Omega}$ are ascribed. The strong form for the thermal problem then reads as:  Find $T: \Omega  \rightarrow \mathbb{R}$ such that $\forall$ $m\in \mathcal{M}$
\begin{equation}
\begin{split}
\label{eq:MMpoisson-strong}
    -\grad \bm{\cdot} ( \kappa^m (\bm{x})\grad T ) = f   ~~~ & \text{ in } \Omega^m \text{,} \\
    [\![ T]\!]^{km} = 0   ~~~ &  \text{ on all } \Gamma_{km}, \\
    [\![ \bm{q}]\!]^{km} = 0   ~~~ &  \text{ on all } \Gamma_{km}, \\ 
     T  = \overline{T}  ~~~ & \text{ on } \Gamma_{\overline{T}}^m \text{,}\\
\end{split}
\end{equation}
where $\Gamma_{km }= \overline{\Omega}^k \cap \overline{\Omega}^m \neq \emptyset$, with $k\in\mathcal{M}$ and $ k \neq m$, are the material interfaces, and $[\![ \cdot]\!]^{km} = (\cdot)^{k} - (\cdot)^{m}$ is the jump of a given quantity over an interface $\Gamma_{km}$.  The material fields are defined $T^m =T(\bm{x})$, $\bm{x} \in \Omega^m$, and  $ \bm{q}^m =  - \kappa^m \grad T^m$. 
The domain $\Gamma_{\overline{T}}^m = \Gamma_{\overline{T}} \cup \partial \overline{\Omega}^m$ are the intersections of the domain boundaries with the material subdomain boundaries. 
$\bm{n}$ denotes the surface normal. 

For this example the method of manufactured solutions is employed using the exact solution 
\begin{align}
    T_{ex}(\bm{x}) = \dfrac{1}{\kappa(\bm{x})} \sin{\dfrac{ 5 \pi (x - 0.2)}{3} }\sin{\dfrac{ 5 \pi y}{3} },
\end{align}
which is constructed to be continuous over the domain, but weakly discontinuous at material interfaces. $T_{ex}$ and its x-derivative are shown in 
Figure \ref{fig:threeMat}. $T_{ex}$ is ascribed as Dirichlet boundary data on the domain boundary. The source term is 
\begin{align}
    f = -\grad \bm{\cdot} ( \kappa(\bm{x})\grad T_{ex}) = \dfrac{25 \pi^2 }{9} \sin{\dfrac{ 5 \pi (x - 0.2)}{3} }\sin{\dfrac{ 5 \pi y}{3} }. 
\end{align}

The discrete form can be defined as: Find $ T^h \in \mathcal{V}^h$ such that $\forall\theta^h\in\mathcal{V}^h$,
\begin{align}
 \sum_{m=1}^n \left[\int_{\Omega^m} \! \! \kappa \grad T^h \bm{\cdot} \grad \theta^h d\Omega\right]   - \int_{\Omega} f \theta^h d\Omega  = \mathcal{R}^{D}\text{  ,}\label{eq:heat-disc-strong}
\end{align}
where $\mathcal{R}^{D}$ is the Nitsche's method residual for enforcing Dirichlet boundary conditions, 
\begin{align} \label{eq:temp-nitsches}
    \nonumber \mathcal{R}^{D} =   \sum_{m=1}^n \Bigg [\mp \int_{\Gamma_{\overline{T}}^m} \kappa (T^h - \overline{T}) (\grad \theta^h \cdot\bm{n}) \,d\Gamma \quad - \\
    \int_{\Gamma_{\overline{T}}^m} \kappa \theta^h(\grad T ^h\cdot\bm{n})\,d\Gamma \quad + \int_{\Gamma_{\overline{T}}^m}\frac{\beta^{D} \kappa }{h}(T^h - \overline{T} ) \theta^h\,d\Gamma \Bigg], 
\end{align}
$h$ is taken as the foreground mesh size, and $\beta^D$ is a user specified constant. For this example $\beta^D = 10$. Once again the symmetric variant of Nitsche's method is employed such that the first quantity in the summation is negative. 
Heaviside enrichment results in strongly discontinuous basis functions, thus continuity of the solution field must be enforced at the interfaces.This is done through the addition of residual $\mathcal{R}^{I}$ to equation \ref{eq:temp-nitsches}. The residual is constructed as 
\begin{align}
    \nonumber \mathcal{R}^{I} =  \sum\limits_{i=1}^n &\sum\limits_{j=i+1}^n \Bigg[- \int_{\Gamma_{ij}} [\![ T^h ]\!] \{ \kappa \grad \theta^h \} \cdot\bm{n}) \,d\Gamma \\
    \nonumber &\quad - \int_{\Gamma_{ij}} [\![ \theta^h ]\!] \{\kappa \grad T ^h\} \cdot\bm{n})\,d\Gamma \\
    &\quad + \int_{\Gamma_{ij}} \gamma^{ij}_T  [\![ T^h ]\!] [\![ \theta^h]\!] \,d\Gamma\Bigg],
\end{align}
where $\{\cdot\} = w^{i}(\cdot)^{i} - w^{j}(\cdot)^{j}$ is the weighted average. 

The function space used to discretize the temperature field is constructed of RKPM shape functions. The shape functions, of polynomial order $n=1$ or $n=2$, have circular supports of normalized radius $a = n+1$. The support radius is normalized by each domain kernel spacing, computed using the distance between the node and its fourth farthest neighbor. Cubic splines are used as kernel functions $\phi$. Interpolation is performed with a boundary conforming quadrilateral mesh, with element sizes equal to nodal spacing and equal order foreground polynomials ($k = n$). 

For classic RKPM, integration is performed using Gauss quadrature points defined upon a uniform grid or equal resolution to the RKPM-nodal spacing. A $6\times 6$ Gauss point rule was employed for linear basis functions and an $8\times 8$ Gauss point rule was employed for quadratic basis functions. For the enriched classic RKPM analysis, the exact characteristic function for each material interface was used from the domain description in equation \ref{eq:analyticDomain}. 

\begin{figure}[t]
	\centering
	\includegraphics[width= \textwidth]{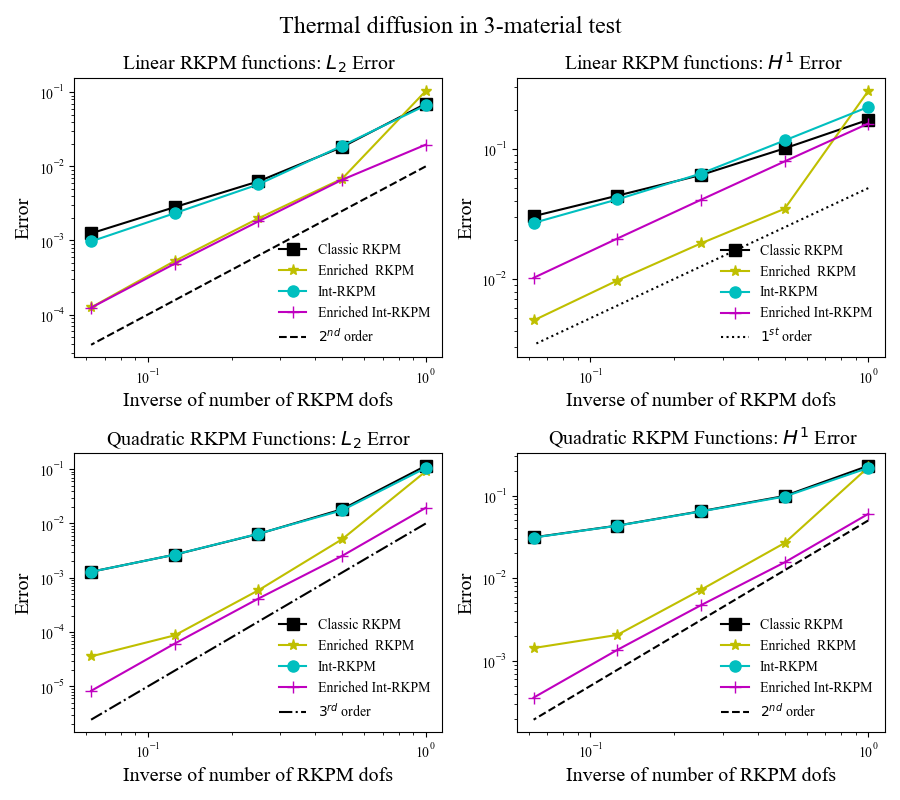}\
	\caption{\label{fig:3matPlot}Convergence data for the three material heat conduction problem, comparing .}
\end{figure}

\begin{figure}[t]
	\centering
	\includegraphics[width= \textwidth]{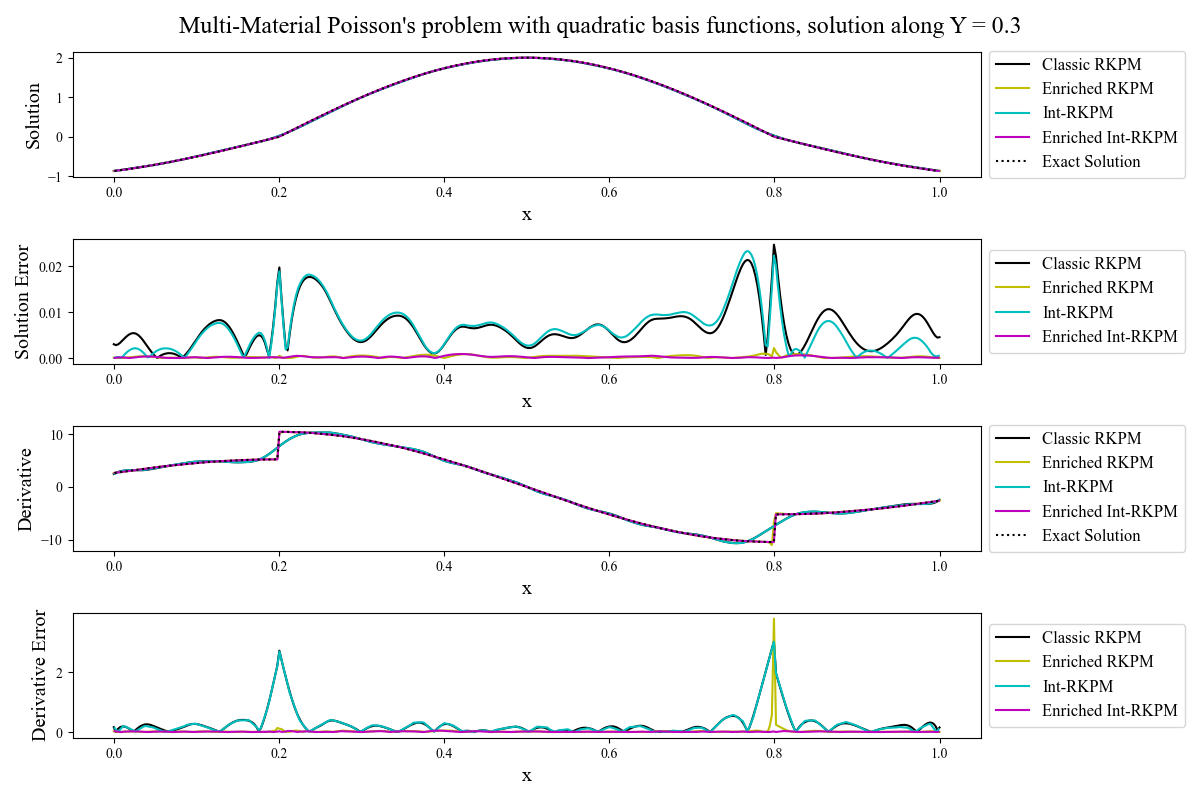}\
	\caption{\label{fig:midlineplots} Enriched RKPM and classic RKPM are compared with their interpolation counterparts. Enriched RKPM and Enriched Int-RKPM accurately resolves the weakly discontinuous temperature field.}
\end{figure}

The results of this convergence study are shown in Figure \ref{fig:3matPlot}. Classic RKPM and Int-RKPM behaved almost identically in the presence of solution discontinuities, see the errors plotted along the line $y=0.3$ in Figure \ref{fig:midlineplots}. Without enrichment, suboptimal error convergence rates are observed due to solution errors near the discontinuities. 

With Enriched RKPM and Enriched-Int RKPM, ideal error convergence rates are achieved for $n=1$, with Enriched-Int RKPM having slightly worse error magnitude in the $H^1$ error norm, following what was observed in subsection \ref{subsec:Poisson}. For $n=2$, Enriched RKPM demonstrates ideal error convergence rates in the pre-asymptotic regime, but rates increase with refinement due to slight errors near the discontinuities (see the yellow line in Figure \ref{fig:midlineplots}).  Enriched Int-RKPM provides a smoother solution space and produces ideal error convergence rates in the asymptotic regime. 

\subsection{Material interfaces with irregular geometries }\label{subsec:eigenstrain}

Employing both the enrichment strategy introduced in Section \ref{subsec:3mat} and the local refinement techniques from Section \ref{subsec:holeInPlate}, Int-RKPM is extended to multi-material PDEs with irregularly shaped interfaces. This capability is verified with a classic benchmarking problem: a circular inclusion is embedded in a semi-infinite medium and a uniform isotropic eigenstrain is imposed upon the inclusion. 

For this example, the inclusion is comprised of Material $1$ with Lam\'{e} constants $\lambda_1= 497.16$ and $\mu_1= 390.63$, while the exterior plate is made of Material $2$ with $\lambda_2 = 656.79$, and $\mu_2 = 338.35$.  The eigenstrain $\varepsilon_0=0.1$ is imposed on the inclusion. The same geometric discretizations and RKPM functions used in Section \ref{subsec:holeInPlate} are employed here, such that $R=1$. To approximate an infinite domain, the exact solution is imposed as Dirichlet boundary data on the top and right of the domain, and symmetry conditions are enforced on the left and bottom.

The strong and discrete forms of the multi-material linear elasticity PDE are omitted from this work for the sake of brevity. The forms are given in full in Appendix B of the author's previous work \cite{fromm_interpolation-based_2024}, and follow a similar form as the multi-material temperature diffusion PDE in section \ref{subsec:3mat}. For this example the mechanical strain is computed by the principle of superposition,
\begin{equation}
    \bm{\varepsilon}_m (\bm{u}) = \begin{cases}
        \bm{\varepsilon}_u(\bm{u}) - \varepsilon_0 \bm{I} \text{, } &\bm{x} \in \Omega_1\\
        \bm{\varepsilon}_u(\bm{u}) \text{, } &\bm{x} \in \Omega_2\\\end{cases}
\end{equation} 
where the total strain $\bm{\varepsilon}_u(\bm{u}) = \frac{1}{2} ( \grad \bm{u} + (\grad \bm{u})^\text{T})$. 

The weakly discontinuous analytic solution for the radial displacement as given in \cite{wang_homogenization_2003} is 
\begin{equation}
    u_r = \begin{cases}
    C_1 r, & r\leq R, \\
    C_1 \dfrac{R^2}{r}, & r\geq R, \\
    \end{cases}
\end{equation}
where $R$ is the radius of the inclusion and 
\begin{equation}
    C_1 = \dfrac{(\lambda_1 + \mu_1)\varepsilon_0}{\lambda_1 + \mu_1 + \mu_2}.
\end{equation}

With enrichment and double interpolation, the midground mesh must also be enriched to accommodate for discontinuities within basis functions. As with the background RKPM functions, midground shape functions with support straddling the material interface are identified.  By restricting the midground function space to discontinuous Galerkin type spaces, each cell is associated with a unique set of degrees of freedom. Thus all cells on the midground mesh $\omega_{e}^{mg}$, covered by triangular elements on the foreground mesh will support basis functions requiring enrichment. 

\begin{figure}[t]
	\centering
	\includegraphics[width= \textwidth]{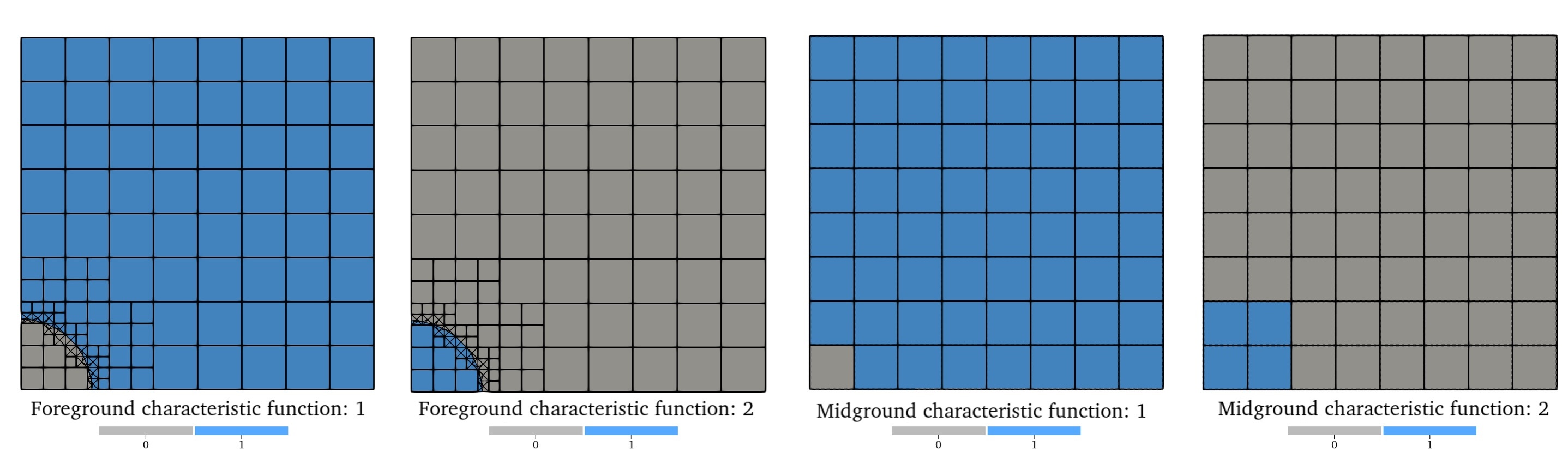}\
	\caption{\label{fig:charFuncs} The material subdomain indicator functions discretized on both the foreground and midground meshes. Foreground meshes use two levels of local refinement.}
\end{figure}

Using the characteristic functions given in Equation \ref{eq:charFunc}, the midground functions are enriched such that 
\begin{align}
    N^{mg, m }_{i} = \psi^m(\bm{x}) N^{mg}_i, \qquad \forall m \in \{1, ..., L\} . 
\end{align}
As this will add additional degrees of freedom to the midground basis, an additional set of characteristic functions must be defined to facilitate enrichment of the background RKPM basis:
\begin{equation}
    \psi^m_{mg} (\bm{x}) = \begin{cases}
        1, \text{ if } \bm{x} \in \{ \omega_{e}^{mg}\}_{m}\\
        0, \text{ else, }
    \end{cases}
\end{equation}
where $ \{ \omega_{e}^{mg}\}_{m}$ is the set of all elements in the midground mesh covered by quadrilateral elements in $\Omega^m$ on the foreground mesh or by triangular elements on the foreground mesh. The foreground and midground characteristic functions are illustrated for the circular inclusion problem in Figure \ref{fig:charFuncs}.

Enriched RKPM functions evaluated for double interpolation are thus expressed as 
\begin{align}\label{eq:enrichB}
    \Psi^m_I(\bm{x}) = \psi_{mg}^m (\bm{x})\Psi_I (\bm{x}), \ \ \forall \ m \in \{1, ..., L\}.
\end{align} 
Note that as the midground material indicator functions are not boundary conforming and may overlap one another, the summation of enriched RKPM functions is no longer equivalent to the original function.

\begin{figure}[t]
	\centering
	\includegraphics[width= \textwidth]{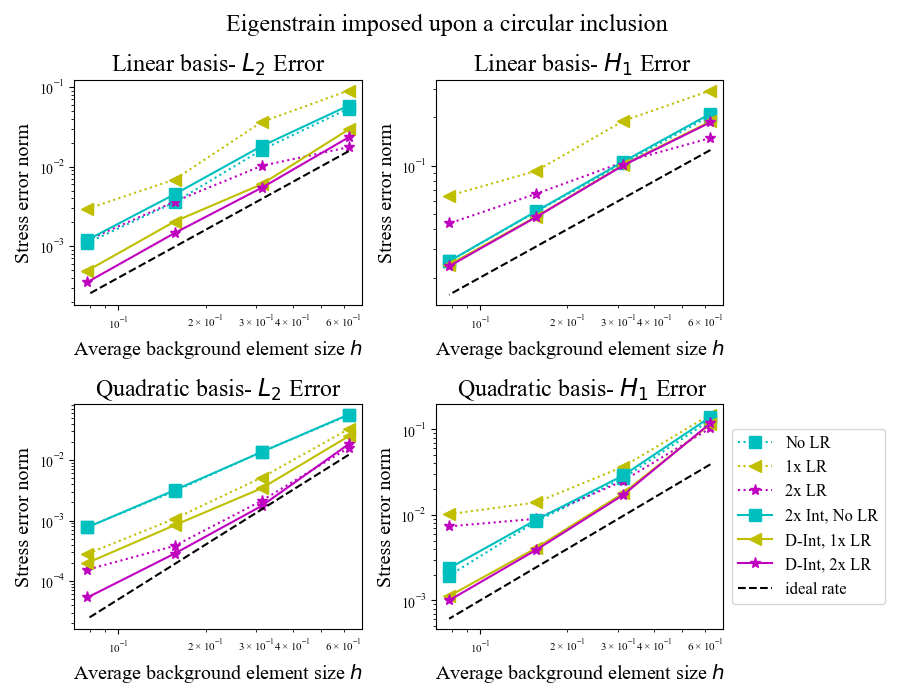}\
	\caption{\label{fig:eigenstrainPlot} Convergence data for eigenstrain imposed upon a circular inclusion, comparing both local refinement and the classic interpolation and double interpolation techniques.}
\end{figure}

With enrichment and double interpolation, ideal error convergence rates are observed for both linear and quadratic functions as shown in Figure \ref{fig:eigenstrainPlot}. Local refinement is required to reduce geometric error to below approximation error for quadratic functions. This problem validates the use of Int-RKPM for multi-material PDEs on irregular geometries. 


\section{Summary }\label{sec:conclusion}

This work introduces a novel interpolation-based reproducing kernel particle method (Int-RKPM). This method is implemented through the opensource finite element software FEniCSx, and all code is publicly available at \cite{fromm_Int-RKPM_2025}. Classical RKPM shape functions are interpolated as linear combinations of Lagrange polynomial shape functions. A foreground mesh is required for integration purposes, but the mesh is not subject to the typical quality constraints of classic finite element analysis. This foreground mesh is used to create a foreground Lagrange polynomial basis for interpolation of the RKPM basis. The interpolated basis retains the vital properties of the original basis, and is suitable for the approximation of PDEs, as illustrated in the numerical examples. 


Unlike the standard RKPM that requires special domain integration techniques such as the stabilized conforming nodal integration \cite{chen_stabilized_2001,chen_non-linear_2002} and variationally consistent integration \cite{hillman_accelerated_2016}  to achieve optimal convergence, the proposed interpolation-based RKPM can be integrated by the standard Gauss integration for optimal convergence, significantly simplified the domain integration complexity. The optimal convergence in Int-RKPM is achieved for reproducing kernel with both linear and quadratic bases, and for solving second-order and forth-order PDEs, demonstrating the robustness of the approach. The geometric topological complexity can be addressed with standard meshing technique without restrictions on mesh quality and hanging nodes for Int-RKPM as they are used for integration purposes. Further, the foreground h- and p- refinements for reducing the geometry and interpolation errors do not affect the number of RKPM approximation coefficients, thus adds no additional cost in solving the Int-RKPM discrete linear system.  These properties allow Int-RKPM to be implemented with standard FEM software. 

Future work will extend the multi-material capabilities of Int-RKPM, in particular investigating the potential of kernel scaling \cite{susuki_image-based_2024,wang_support_2024} to represent discontinuities as opposed to enrichment. Kernel scaling has the potential to represent either weak or strong discontinuities, avoiding the necessities of Nitsches interface conditions.

\section*{Acknowledgments}
The authors were supported by National Science Foundation award number 2104106.


\bibliographystyle{unsrt} 
\bibliography{int-based-RKPM} 
\end{document}